\documentclass{iopart}
\usepackage{iopams,amsbsy,eufrak,mathrsfs,graphicx}  

\def\espaitemps{({\cal V},g)}
\def\varietat{{\cal V}}

\def\S{\Sigma}

\def\be{\begin{equation}}
\def\ee{\end{equation}}
\def\bea{\begin{eqnarray}}
\def\eea{\end{eqnarray}}
\def\bean{\begin{eqnarray*}}
\def\eean{\end{eqnarray*}}

\def\proof{\noindent{\em Proof.\/}\hspace{3mm}}
\def\fin{\hfill \rule{2.5mm}{2.5mm}\\ \vspace{0mm}}
\def\finn{\hfill \rule{2.5mm}{2.5mm}}

\newtheorem{prop}{Proposition}[section]
\newtheorem{theorem}{Theorem}[section]

\newtheorem{lemma}{Lemma}[section]
\newtheorem{coro}{Corollary}[section]
\newtheorem{conj}{Conjecture}

\begin{document}

\title{2nd-order symmetric Lorentzian manifolds I: characterization
  and general results}
\author{Jos\'e M.M. Senovilla}
\address{F\'{\i}sica Te\'orica, Universidad del Pa\'{\i}s Vasco,
Apartado 644, 48080 Bilbao, Spain \\

E-mail: josemm.senovilla@ehu.es}
\begin{abstract} 
The $n$-dimensional Lorentzian manifolds with vanishing second 
covariant derivative of the Riemann tensor --- 2--symmetric 
spacetimes--- are characterized and classified. The main result is
that either they are locally symmetric or they
have a covariantly constant null vector field, in this case defining a 
subfamily of Brinkmann's class in $n$ dimensions.
Related issues and applications
are considered, and new open questions presented. 
\end{abstract} 

Keywords: Symmetric spaces, semi-Riemannian manifolds, parallel null 
vector fields, curvature invariants, General Relativity.

MSC: 53B30, 53B20, 53C50, 53Z05

\vspace{1cm}
\section{Introduction}
\label{intro}
The aim of this paper is to characterize the 
manifolds $\varietat$ with a metric $g$ of Lorentzian signature such 
that the Riemann tensor $R^{\alpha}{}_{\beta\gamma\delta}$ of 
$\espaitemps$  {\em locally} satisfies the second order condition 
\be
\nabla_{\mu}\nabla_{\nu} R^{\alpha}{}_{\beta\gamma\delta}=0. 
\label{DDR}
\ee
It is quite surprising that, hitherto, despite their simple 
definition, this type of Lorentzian manifolds have been hardly 
considered in the literature. Probably this is due to a combination
of reasons of diverse origin, such as:
\begin{itemize}
\item the classical Riemannian
result \cite{L,NO,Ta} according to which all such spaces are actually locally
symmetric \cite{C} if the metric is positive definite, see
section \ref{parafernalia}; 
\item even in cases with other signatures for the metric, there are
well-known results \cite{Ta,CS1} restricting the possibility of
(\ref{DDR}). 
Thus, one can say that manifolds with
the property (\ref{DDR}) are somehow exceptional: see next section for
a recollection of such results and some new expanded ones along these lines.
\item characterizing the exceptional cases left over by the previous
considerations turns out to be, if not difficult, very laborious
indeed. Given that one knows beforehand that the result will be a very
special class of metrics this may have prevented some researchers from
considering the problem. It must be remarked, nevertheless, that the
4-dimensional case is easily
solved ---see the footnote number 1 for the sketch of the solution!
\end{itemize}

The interest of these spacetimes comes also from different
perspectives and fields. To start with, they have an obvious
mathematical interest. Having such a simple and natural definition,
they should be identified for arbitrary
signature. Furthermore, there are direct geometric interpretations of condition (\ref{DDR}) in
analogy with the case of a vanishing first covariant derivative \cite{O}: the
tensor field $\nabla_{\nu} R^{\alpha}{}_{\beta\gamma\delta}$ is invariant under parallel
displacement, and therefore one can say that the curvature is locally
a ``linear'' function of appropriate coordinates. This has a relation
to the local holonomy group of the manifold. 

In a more concrete manner, (i) the condition
(\ref{DDR}) implies that given any geodesic curve on the manifold with
tangent vector $\vec v$, and parallelly propagated vector fields $\vec
X,\vec Y,\vec Z$ along the geodesic, then the vector field
$$
\nabla_{\vec v}R(\vec X,\vec Y)\vec Z
$$
is itself parallelly propagated along the geodesic. One could say,
leaving aside rigour, that the vector field $R(\vec X,\vec Y)\vec Z$ is a ``linear''
function on the affine parameter of the geodesics. And similarly,
(ii) take the sectional curvature $\kappa$ of the manifold \cite{O} at
any point relative to the tangent planes at that point. Let
the tangent planes be parallelly propagated along geodesics, then the
function $\nabla_{\vec v}\kappa$ remains constant. 

Nevertheless, there is an important difference between (\ref{DDR}) and
its first order version with only one covariant derivative: in the
latter case, it is known that there is (locally) a geodesic symmetry
\cite{O}, that is to say, there is an isometry $\varphi$ of the
Lorentzian manifold which acts on the tangent spaces as $d\varphi
=-$Id. In the case of (\ref{DDR}), there is nothing of this sort in
general, as one would need an isometry whose {\em square}, when acting
on the tangent spaces, were not proportional to Id. But this is
impossible in general, even in 2-dimensional cases.

Having considered their immediate geometric interpretation, the spaces
satisfying (\ref{DDR}) may also have important applications in physics,
and in particular in theories concerning gravitation whose classical
arena is a Lorentzian manifold. Just to cite a few possible
applications or areas of potential interest, let us mention the
following:
\begin{itemize}
\item they are of interest in the branch of invariant
classification of Lorentzian manifolds, and of the ``equivalence
problem'' \cite{Exact}, that is, to decide if two given spacetimes are
locally the same.
\item in analogy with flat or conformally flat spacetimes, which can
be characterized by the vanishing of the Bel or the Bel-Robinson tensors \cite{S}
respectively, the spacetimes with a vanishing ``(super)$^2$-energy''
have $\nabla_{\nu} R^{\alpha}{}_{\beta\gamma\delta}=0$, 
and those with zero ``(super)$^3$-energy'' will satisfy
(\ref{DDR}). Therefore, these spacetimes can shed some light into the potentially
high number of conserved quantities that one can form with those
tensors, see \cite{S}, and whether the superenergy construction stops
effectively at some level.
\item these spacetimes may arise whenever expansions in normal 
coordinates are used, or relevant. The higher order
terms in the expansions contain summands proportional to the higher order
derivatives of the Riemann tensor. Therefore, if these derivatives happen to vanish the 
expansions may become more manageable. Similarly, if one considers the local form of a Lorentzian manifold around any null geodesic (the so-called ``Penrose limits''
\cite{Plimit}) the various levels of approximation contain summands
correspoding to higher-order conditions of type (\ref{DDR}). It is worth mentioning here that the main result obtained in this paper will actually imply that many of the 
non-locally-symmetric spacetimes satisfying (\ref{DDR}) will
actually be of Penrose limit type: a plane wave
\cite{EK,Plimit,Exact}.
\item it will also turn out that they are closely related to
spacetimes with special properties concerning its curvature
invariants. Many curvature scalar invariants will be zero or constant,
linking these spaces to some families of spacetimes which have
become a focus of recent interest, see \cite{PPCM,CMPP2,CHP} and references
therein.
\item related to the previous point, they can be of relevance as
solutions in higher order Lagrangian physical
theories including gravity, as then only a
finite number of terms are relevant. In particular, they can be examples of exact
solutions (for the background) in string theory when supported by
appropriate matter contents, see e.g.\ \cite{HT,PPCM,Coley} and references therein.
\item actually, they can  provide examples of exact solutions for
backgrounds for 11-dimensional supergravity and relatives via M-theory
\cite{F-O,BFHP}. This is in fact interrelated with the Penrose limits
mentioned previously \cite{Gu,BFHP,BFHP2}.
\item in the previous point as well as in general solitonic or
black-hole solutions of supersymmetric theories, these solutions
must have a covariantly constant spinor, leading to covariantly
constant null vector or tensor fields, see \cite{F-O}. These will be seen to arise
naturally in the spacetimes satisfying (\ref{DDR}).
\item they can provide some examples of solutions in field equations
for which the counter-terms regularizing quantum fluctuations are
vanishing \cite{Gib}.
\item they also seem to be related to what has been called the
``$\epsilon$-property'' \cite{PCM}, which refers basically to the
possibility that some components of the Riemann tensor (and its
derivatives) can be made as small as desired by a judicious choice of
basis. 
\end{itemize}

Observe that in proper Riemannian spaces (i.e., with a 
positive-definite metric), the vanishing of the square of a tensor 
implies that the tensor itself is zero. Though this is not so for the 
semi-Riemannian case, it will be shown that, {\it in the case of 
Lorentzian signature}, the use of causal tensors and superenergy 
techniques, see \cite{S,BS}, provides a valid and productive alternative and similar 
results can be found. I will make use of this in several places.

The main result of this paper is that all Lorentzian manifolds
satisfying (\ref{DDR}) but with a non-vanishing first derivative of
the Riemann tensor must necessarily have a null covariantly constant
vector field. Therefore, all these spacetimes  belong to a general class of
Lorentzian manifolds known as the Brinkmann class \cite{Br}.
As a matter of fact, this result can be obtained rather quickly
in four dimensions by using spinors,\footnote{Indeed, if the
curvature satisfies (\ref{DDR}), using the Ricci identities for the
curvature spinors ---formulae (4.9.13-15) in \cite{PR}--- one easily derives
that the Weyl spinor satisfies the condition appearing at the
beginning of page 261 in \cite{PR}, leading to a type N Weyl
tensor. The vanishing covariant derivative of the unique principal
null direction follows then in various possible ways, for example
using the superenergy tensors of \cite{S}.} but as we
will see it is far from
obvious in higher dimensions. 

Of course, the fundamental 
reason behind this result is the now known and better understood 
degenerately decomposability \cite{Wu,Wu1,Wu2,Hall} of the 
spacetime in this case, which was completely analyzed in \cite{CW} 
for the Lorentzian case, and has been a subject of recent interest 
with many interesting applications is supergravity and string theory, see e.g. 
\cite{F-O,Bry,FP,BFHP,BK} and references therein. This will be briefly
analyzed in subsection \ref{holonomy}.

Throughout, $\varietat$ will denote a differentiable manifold and $g$
a metric tensor of sufficient differentiability. When considering
Lorentzian signature, the choice of signature will be
$(-,+,\dots,+)$, so that timelike vectors have a negative length,
spacelike vectors a positive one, and the null vectors, which include
the zero vector, have vanishing length. Lorentzian manifolds are
assumed to be time orientable with a chosen future direction. Given
the type of formulas needed in the paper I have preferred to write
most of the calculations and expressions using index notations. It
should be clearly understood, however, that no particular basis has
been chosen, and the results are general. They can of course be
rewritten in index-free form if desired.

\section{Generalizations of symmetric spaces. Generic 
cases}
\label{parafernalia}
Semi-Riemannian (or pseudo-Riemannian) manifolds satisfying 
(\ref{DDR}) constitute an obvious generalization of the well-known 
locally {\em symmetric} spaces which satisfy
\be
\nabla_{\mu} R^{\alpha}{}_{\beta\gamma\delta}=0 \label{DR}
\ee
and were introduced, largely studied and classified by \'{E}. Cartan 
\cite{C} in the proper 
Riemannian case , see e.g. \cite{C1,KN,H}, and later in 
\cite{CM,CW,CP,CK} for the Lorentzian case using results from
\cite{Ber}. The general semi-Riemannian 
cases were treated in e.g. \cite{CP1,O}, but these results seem to be
incomplete, see the recent studies in \cite{KO} and references therein. Locally 
symmetric spaces are themselves generalizations of the constant 
curvature spaces and, as a matter of fact, there is a natural 
hierarchy of conditions, shown in Table \ref{ta}, that can be placed 
on the curvature tensor. In that table, the restrictions on the 
curvature tensor decrease towards the right and each class is 
strictly contained in the following ones. Thus, all constant 
curvature manifolds are obviously symmetric but the converse is not 
true, and analogously for the other cases. The table has been stopped 
at the level of semi-symmetric spaces, defined by the condition 
$\nabla_{[\mu}\nabla_{\nu]} R^{\alpha}{}_{\beta\gamma\delta}=0$, 
where round and square brackets enclosing indices indicate 
symmetrization and antisymmetrization, respectively. As a matter of fact, 
semi-symmetric spaces were introduced also by Cartan \cite{C1} and 
studied in \cite{Sz,Sz1} as the natural generalization of symmetric 
spaces for the proper Riemannian case, see also \cite{BKV} and 
references therein;  in the 4-dimensional Lorentzian case they were 
studied in \cite{K}, and also in e.g. \cite{DVV,DDVV,DDSVY,HV} and 
references therein. Explicit proof that there exist non-symmetric 
semi-symmetric spaces was given in \cite{T}.

\begin{table}
\caption{The hierarchy of conditions on the Riemann tensor}
\begin{center}
\begin{tabular}{|l|c|c|l|}
\hline
$R^{\alpha}{}_{\beta\gamma\nu}\propto 
\delta^{\alpha}_{\gamma}g_{\beta\nu}-\delta^{\alpha}_{\nu}g_{\beta\gamma}$ 
& $\nabla_{\mu} R^{\alpha}{}_{\beta\gamma\delta}=0$ &
$\nabla_{\mu}\nabla_{\nu} R^{\alpha}{}_{\beta\gamma\delta}=0$ & 
$\nabla_{[\mu}\nabla_{\nu]} 
R^{\alpha}{}_{\beta\gamma\delta}=0$ \\
\hline
constant curvature & symmetric & 2-symmetric & semi-symmetric \\
\hline
\end{tabular}
\label{ta}
\end{center}
\label{default}
\end{table}%

Two obvious questions arise: 
\begin{enumerate}
\item why is semi-symmetry considered to 
be the natural generalization of local symmetry, instead of 
(\ref{DDR})?
\item why not go on further to higher derivatives of 
the Riemann tensor? 
\end{enumerate}
The answer to both questions is actually the 
same: a classical theorem \cite{L,NO,Ta} states that in any proper 
Riemannian manifold
\be
\nabla_{\mu_1}\dots \nabla_{\mu_k}R^{\alpha}{}_{\beta\gamma\delta}=0 
\hspace{1cm} \Longleftrightarrow  \hspace{1cm} \nabla_{\mu} 
R^{\alpha}{}_{\beta\gamma\delta}=0 \label{DkR}
\ee
for any $k\geq 1$ so that, in particular, (\ref{DDR}) is strictly 
equivalent to (\ref{DR}) in proper Riemannian spaces. 

This may also be the reason of why there seems to be no name 
for the condition (\ref{DDR}) in the literature. However, an 
analogous condition has certainly been used for the so-called 
recurrent spaces: if there exists $A_{\mu_1\dots \mu_k}$ such that 
$\nabla_{\mu_1}\dots\nabla_{\mu_k}R^{\alpha}{}_{\beta\gamma\delta}=A_{\mu_1\dots 
\mu_k}R^{\alpha}{}_{\beta\gamma\delta}$ then the space is called 
$k$-recurrent (e.g. \cite{Tak,CS1}), in particular second order 
recurrent (or 2-recurrent) for $k=2$ \cite{L,R1,Th,Th2,MTh,Tak,CS} 
and recurrent for $k=1$ (see e.g. \cite{RWW} and references therein). 
Thus, I will call the spaces satisfying (\ref{DDR}) {\em second-order 
symmetric}, or in short {\em 2-symmetric}, and more generally 
$k$-symmetric when the left condition in (\ref{DkR}) holds. The whole
class of $k$-symmetric spacetimes for all $k>1$ has been called ``higher
order symmetric spaces'' recently in \cite{PCM}.

\subsection{Results at {\em generic} points}
\label{generic}
As a matter of fact, the equivalence (\ref{DkR}) holds as well in 
``generic" cases of semi-Riemannian manifolds of any signature. For 
some results on this one can consult \cite{Ta,CS1}. A typical 
reasoning would be as follows. Assume that the left side of 
(\ref{DkR}) holds for either $k=2,3$, then the Ricci identity applied 
to 
$\nabla_{[\lambda}\nabla_{\mu]}\nabla_{\nu}R_{\alpha\beta\gamma\delta}$ 
provides
\be
R^{\rho}{}_{\nu\lambda\mu}\nabla_{\rho}R_{\alpha\beta\gamma\delta}+
R^{\rho}{}_{\alpha\lambda\mu}\nabla_{\nu}R_{\rho\beta\gamma\delta}+
R^{\rho}{}_{\beta\lambda\mu}\nabla_{\nu}R_{\alpha\rho\gamma\delta}+
R^{\rho}{}_{\gamma\lambda\mu}\nabla_{\nu}R_{\alpha\beta\rho\delta}+
R^{\rho}{}_{\delta\lambda\mu}\nabla_{\nu}R_{\alpha\beta\gamma\rho}=0 
\label{RicDR}
\ee
so that if at any point $p\in \varietat$ the matrix 
$(R^{\alpha\beta}{}_{\gamma\delta})|_p$ of the Riemann tensor, 
considered as an endomorphism on the space of 2-forms $\Lambda_2(p)$, 
is non-singular, then we can multiply (\ref{RicDR}) by the inverse 
matrix of $(R^{\alpha\beta}{}_{\gamma\delta})|_p$ getting at $p$
$$
\delta^{\nu}_{[\lambda}\nabla_{\mu]}R^{\alpha\beta\gamma\delta}+
\delta^{\alpha}_{[\lambda}\nabla^{\nu}R_{\mu]}{}^{\beta\gamma\delta}+
\delta^{\beta}_{[\lambda}\nabla^{\nu}R^{\alpha}{}_{\mu]}{}^{\gamma\delta}+
\delta^{\gamma}_{[\lambda}\nabla^{\nu}R^{\alpha\beta}{}_{\mu]}{}^{\delta}+
\delta^{\delta}_{[\lambda}\nabla^{\nu}R^{\alpha\beta\gamma}{}_{\mu]}=0
$$
which after contracting $\nu$ with $\lambda$ leads, using the second 
Bianchi identity, to
\be
(n+1)\nabla_{\mu}R_{\alpha\beta\gamma\delta}=
g_{\alpha\mu}\nabla_{\rho}R^{\rho}{}_{\beta\gamma\delta}+
g_{\beta\mu}\nabla_{\rho}R_{\alpha}{}^{\rho}{}_{\gamma\delta}+
g_{\gamma\mu}\nabla_{\rho}R_{\alpha\beta}{}^{\rho}{}_{\delta}+
g_{\delta\mu}\nabla_{\rho}R_{\alpha\beta\gamma}{}^{\rho}\, .\label{pas}
\ee
Contracting here with $g^{\mu\alpha}$ we easily get 
$\nabla_{\rho}R^{\rho}{}_{\beta\gamma\delta}=0$ and introducing this 
result into (\ref{pas}) we finally arrive at 
$\nabla_{\mu}R_{\alpha\beta\gamma\delta}=0$. Thus, 2-symmetry (or 
3-symmetry) implies local symmetry on a neighbourhood of any $p\in 
\varietat$ at which the Riemann tensor matrix 
$(R^{\alpha\beta}{}_{\gamma\delta})|_p$ is non-singular. This is the 
meaning of the word generic used above.

Keeping the meaning of the word ``generic'' in mind, this type of reasoning 
can be extended to arbitrary tensor fields---with possible stronger 
results depending on their order and symmetries---. For instance, one can 
prove the following general result.
\begin{prop}
\label{2implies1}
Let $\tilde T$ be any tensor field. Around \underline{{\em generic}} points, the vanishing of 
its second covariant derivative implies the vanishing 
of its first covariant derivative.
\end{prop}
\proof Let us denote by $T_{\alpha_1\dots \alpha_q}$ the totally 
covariant tensor field equivalent to $\tilde T$ by lowering all 
contravariant indices. Assume 
that $\nabla_{\lambda}\nabla_{\mu}T_{\alpha_1\dots \alpha_q}=0$, then the 
Ricci identity applied to 
$\nabla_{[\lambda}\nabla_{\mu]}T_{\alpha_1\dots \alpha_q}=0$ and to 
$\nabla_{[\lambda}\nabla_{\mu]}\nabla_{\nu}T_{\alpha_1\dots 
\alpha_q}=0$ provides, respectively
\bea
\sum_{i=1}^{q} R^{\rho}{}_{\alpha_i\lambda\mu}
T_{\alpha_1\dots\alpha_{i-1}\rho\alpha_{i+1}\dots\alpha_q}=0, 
\label{RicT}\\
R^{\rho}{}_{\nu\lambda\mu}\nabla_{\rho}T_{\alpha_1\dots 
\alpha_q}+\sum_{i=1}^{q} R^{\rho}{}_{\alpha_i\lambda\mu}\nabla_{\nu}
T_{\alpha_1\dots\alpha_{i-1}\rho\alpha_{i+1}\dots\alpha_q}=0.\label{RicDT}
\eea
At any $p$ where $\det (R^{\alpha\beta}{}_{\gamma\delta})\neq 0$ then 
we have
\bean
\sum_{i=1}^{q}\left(g_{\alpha_i\lambda}
T_{\alpha_1\dots\alpha_{i-1}\mu\alpha_{i+1}\dots\alpha_q}-
g_{\alpha_i\mu}T_{\alpha_1\dots\alpha_{i-1}\lambda\alpha_{i+1}\dots\alpha_q}\right)=0,
\hspace{1cm}\\
g_{\nu\lambda}\nabla_{\mu}T_{\alpha_1\dots \alpha_q}-
g_{\nu\mu}\nabla_{\lambda}T_{\alpha_1\dots \alpha_q}+
\sum_{i=1}^{q}\left(g_{\alpha_i\lambda}
\nabla_{\nu}T_{\alpha_1\dots\alpha_{i-1}\mu\alpha_{i+1}\dots\alpha_q}-
g_{\alpha_i\mu}\nabla_{\nu}
T_{\alpha_1\dots\alpha_{i-1}\lambda\alpha_{i+1}\dots\alpha_q}\right)=0
\eean
so that covariantly differentiating the first and substracting the 
second one gets
$$
g_{\nu\lambda}\nabla_{\mu}T_{\alpha_1\dots \alpha_q}-
g_{\nu\mu}\nabla_{\lambda}T_{\alpha_1\dots \alpha_q}=0
$$
and contracting here $\nu$ and $\lambda$ one finally proves 
$\nabla_{\mu}T_{\alpha_1\dots \alpha_q}=0$.\finn
\begin{coro}
For any tensor field $T$, and at \underline{{\em generic}} 
points, one has
$$
\stackrel{k}{\overbrace{\nabla\cdots\cdots\nabla}}\,  T =0 \Longleftrightarrow \nabla T =0
$$
for any $k\geq 1$.\fin
\end{coro}

These results apply in particular to the Riemann tensor and, actually,
stronger results can be proven sometimes. As an interesting example,
let us mention that by application of the previous results one can prove a 
conjecture in \cite{CS1}, namely, that all $k$-symmetric (and also 
all $k$-recurrent) spaces are necessarily of constant curvature on a 
neighbourhood of any $p\in \varietat$ at which the Riemann tensor 
matrix is non-singular.
As a matter of fact, let us prove the following slightly more general result
\begin{theorem}
\label{semithenconstant}
All {\em semi-symmetric} spaces are of constant curvature at generic 
points.
\end{theorem}
\proof Assume that $\nabla_{[\lambda}\nabla_{\mu]} 
R^{\alpha}{}_{\beta\gamma\delta}=0$ and apply here the Ricci identity 
to get
\be
R^{\rho}{}_{\alpha\lambda\mu}R_{\rho\beta\gamma\delta}+
R^{\rho}{}_{\beta\lambda\mu}R_{\alpha\rho\gamma\delta}+
R^{\rho}{}_{\gamma\lambda\mu}R_{\alpha\beta\rho\delta}+
R^{\rho}{}_{\delta\lambda\mu}R_{\alpha\beta\gamma\rho}=0 \label{RicR}
\ee
so that as before, at any $p$ with a non-singular Riemann tensor 
matrix,
$$
\delta^{\alpha}_{[\lambda}R_{\mu]}{}^{\beta\gamma\delta}+
\delta^{\beta}_{[\lambda}R^{\alpha}{}_{\mu]}{}^{\gamma\delta}+
\delta^{\gamma}_{[\lambda}R^{\alpha\beta}{}_{\mu]}{}^{\delta}+
\delta^{\delta}_{[\lambda}R^{\alpha\beta\gamma}{}_{\mu]}=0
$$
and contracting here $\alpha$ with $\lambda$ one derives
$(n-1) R_{\mu\beta\gamma\delta}=2g_{\mu[\gamma}R_{\delta]\beta}$ 
immediately implying 
$g_{\mu[\gamma}R_{\delta]\beta}=R_{\mu[\gamma}g_{\delta]\beta}$ from 
where contraction of $\mu$ and $\gamma$ gives 
$(n-1)R_{\beta\gamma}=Rg_{\beta\gamma}$ which together with the 
previous formula for the Riemann tensor proves the result.\fin

In the previous proof, and in the rest of the paper, 
$R_{\mu\nu}\equiv R^{\rho}{}_{\mu\rho\nu}$ and $R\equiv g^{\mu\nu}R_{\mu\nu}$
denote the Ricci tensor and the scalar curvature, 
respectively.

Therefore, there is little room for spaces (necessarily of 
non-Euclidean signature) which are $k$-symmetric but {\em not} 
symmetric nor of constant curvature. It is remarkable that there have 
been many studies on 2-recurrent (or conformally 2-recurrent, see 
subsection \ref{conformal+ricci}) spaces, which certainly include the 2-symmetric ones, 
but surprisingly enough the assumption that they were {\em not} 
2-symmetric was always, either implicitly or explicitly 
\cite{R1,Th,Th1,Th2,MTh,E-K,E-K1}, made. Thus, the present paper tries to 
fill in this gap: the main purpose is to find these special manifolds 
for the case of 2-symmetry and Lorentzian signature. This will be done
from section \ref{lor-bis} onwards. Before that, let us
collect some general results valid in arbitrary signature and some
possible nomenclature concerning the conformal and Ricci curvature tensors.

\subsection{Identities in 2-symmetric semi-Riemannian manifolds}
\label{identities}
Some basic tensor calculation is needed, mainly to prove 
quadratic identities which hold in general 2-symmetric spaces of any
signature and that, with the help of the lemmas in section \ref{lor-bis}, can 
be used to get the sought results. 

To start with, we need  a 
generalization to non-generic points ---i.e. to the case with a possibly degenerate 
Riemann tensor matrix--- of the calculations and results presented in subsection 
\ref{generic}.
\begin{lemma}
Let $\espaitemps$ be an $n$-dimensional {\em 2-symmetric} 
semi-Riemannian manifold of any signature. If 
$\nabla_{\lambda}\nabla_{\mu}T_{\mu_1\dots \mu_q}=0$ then
\bea
\sum_{i=1}^{q} \nabla_{\nu}R^{\rho}{}_{\alpha_i\lambda\mu}
T_{\alpha_1\dots\alpha_{i-1}\rho\alpha_{i+1}\dots\alpha_q}
-R^{\rho}{}_{\nu\lambda\mu}\nabla_{\rho}T_{\alpha_1\dots 
\alpha_q}=0,\label{eq4}\\
(\nabla_{\nu}R^{\rho}{}_{\tau\lambda\mu}+\nabla_{\tau}R^{\rho}{}_{\nu\lambda\mu})
\nabla_{\rho}T_{\mu_1\dots \mu_q}=0, \hspace{2cm}\label{basic}\\
(\nabla_{\nu}R^{\rho}{}_{\mu}-\nabla_{\mu}R^{\rho}{}_{\nu})
\nabla_{\rho}T_{\mu_1\dots \mu_q}=0, \,\,\, 
(\nabla^{\rho}R_{\mu\nu}-2\nabla_{\nu}R^{\rho}{}_{\mu})
\nabla_{\rho}T_{\mu_1\dots \mu_q}=0. \label{treq4}
\eea
\end{lemma}
\proof As (\ref{RicT}) and (\ref{RicDT}) hold, by covariantly 
differentiating the first and substracting the second one gets 
(\ref{eq4}). Computing the covariant derivative of this equation and 
that of (\ref{RicDT}), using the 2-symmetry, and substracting them 
one arrives immediately at (\ref{basic}).
Contracting here $\tau,\nu$, or $\tau,\lambda$, one gets the two in 
(\ref{treq4}), respectively.\finn

Recall the decomposition of the Riemann tensor,
\be
R_{\alpha\beta\lambda\mu}=C_{\alpha\beta\lambda\mu}+\frac{2}{n-2}\left(
R_{\alpha[\lambda}g_{\mu]\beta}-R_{\beta[\lambda}g_{\mu]\alpha}\right)-
\frac{R}{(n-1)(n-2)}\left(g_{\alpha\lambda}g_{\beta\mu}-
g_{\alpha\mu}g_{\beta\lambda}\right)\label{weyl}
\ee
where $C_{\alpha\beta\lambda\mu}$ is the conformal or Weyl curvature
tensor, which satisfies the same symmetry properties as the Riemann
tensor and is also traceless
$$
C_{\alpha\beta\lambda\mu}=C_{[\alpha\beta][\lambda\mu]}, \,\,\,\,
C_{\alpha[\beta\lambda\mu]}=0, \,\,\,\,
C^{\rho}{}_{\beta\rho\mu}=0
$$
where of course the first two imply that
$C_{\alpha\beta\lambda\mu}=C_{\lambda\mu\alpha\beta}$. Then the
following general identities hold for arbitrary semi-symmetric spaces
\begin{lemma}
The Riemann, Ricci and Weyl tensors of an $n$-dimensional {\em 
semi-symmetric} semi-Riemannian manifold of any signature satisfy 
(\ref{RicR}) as well as
\bea
R_{\rho(\mu}R^{\rho}{}_{\nu)\alpha\beta}=0, \,\,\, 
R^{\rho}{}_{\mu[\alpha\beta}R_{\gamma]\rho}=0, \,\,\, 
C^{\rho}{}_{\mu[\alpha\beta}R_{\gamma]\rho}=0, \,\,\, 
R^{\rho\sigma}R_{\rho\mu\sigma\nu}=R_{\mu}{}^{\rho}R_{\rho\nu}, 
\label{RR}\\
R^{\rho}{}_{\alpha\lambda\mu}C_{\rho\beta\gamma\delta}+
R^{\rho}{}_{\beta\lambda\mu}C_{\alpha\rho\gamma\delta}+
R^{\rho}{}_{\gamma\lambda\mu}C_{\alpha\beta\rho\delta}+
R^{\rho}{}_{\delta\lambda\mu}C_{\alpha\beta\gamma\rho}=0, 
\label{RicC}\\
(n-2)\left(C_{\rho[\alpha}{}^{\lambda\mu}C^{\rho}{}_{\beta]\gamma\delta}+
C_{\rho[\gamma}{}^{\lambda\mu}C^{\rho}{}_{\delta]\alpha\beta}\right)-
2\left(R_{[\alpha}{}^{[\lambda}C^{\mu]}{}_{\beta]\gamma\delta}+
R_{[\gamma}{}^{[\lambda}C^{\mu]}{}_{\delta]\alpha\beta}\right)-\nonumber\\
-2\left(R_{\rho}{}^{[\lambda}\delta^{\mu]}_{[\alpha}C^{\rho}{}_{\beta]\gamma\delta}+
R_{\rho}{}^{[\lambda}\delta^{\mu]}_{[\gamma}C^{\rho}{}_{\delta]\alpha\beta}\right)+
2\frac{R}{n-1}\left(\delta^{[\lambda}_{[\alpha}C^{\mu]}{}_{\beta]\gamma\delta}+
\delta^{[\lambda}_{[\gamma}C^{\mu]}{}_{\delta]\alpha\beta}\right)=0 
\label{CC}
\eea
and their non-written traces, such as the appropriate specializations 
of (\ref{treq4}).
\end{lemma}
\proof
The first in (\ref{RR}) is the trace of 
(\ref{RicR}) (or the Ricci identity for the Ricci tensor), the fourth 
is its trace, the second follows from the first by taking a cyclic 
permutation of indices, and the third follows form the second by 
using (\ref{weyl}). On the other hand, (\ref{RicC}) is the Ricci identity for 
the Weyl tensor and finally (\ref{CC}) follows from (\ref{RicC}) on using 
(\ref{weyl}) again.\finn

In the particular case of 2-symmetric manifolds we also have that 
\begin{lemma}
The Riemann, Ricci and Weyl tensors of an $n$-dimensional {\em 
2-symmetric} semi-Riemannian manifold of any signature satisfy 
(\ref{RicDR}) as well as
\be
\nabla_{(\tau}R^{\rho}{}_{\nu)\lambda\mu}\nabla_{\rho}R_{\alpha\beta\gamma\delta}=0, 
\,\,
\nabla_{(\tau}R^{\rho}{}_{\nu)\lambda\mu}\nabla_{\rho}C_{\alpha\beta\gamma\delta}=0, 
\,\,
\nabla_{(\tau}R^{\rho}{}_{\nu)\lambda\mu}\nabla_{\rho}R_{\alpha\beta}=0, 
\label{basicR}
\ee
and their non-written traces, such as the appropriate specializations 
of (\ref{treq4}).
\end{lemma}
\proof
These are the particularization of (\ref{basic}) to the various
curvature tensors.\finn

\subsection{Special vector fields}
\label{homothetic}
We briefly collect an important result, and its consequences,
that has been used classically (see e.g. \cite{Ta,NO})
to study the second order symmetric semi-Riemannian manifolds. This
will be useful later in the proof of Proposition
\ref{constornull}. The idea is to explore the consequences
of the existence of a
vector field whose covariant derivative is proportional to the metric,
so that it has vanishing second covariant derivative.
\begin{lemma}
\label{homo}
Let $\espaitemps$ be 2-symmetric with arbitrary signature.
If a 1-form $v_{\mu}$ satisfies
$$
\nabla_{\nu}v_{\mu}=c\, g_{\mu\nu}
$$
for a constant $c$, then either $c=0$
or the manifold is locally symmetric.
\end{lemma}
\proof The condition of the lemma implies that 
$v^{\mu}$ is a homothetic vector \cite{Yano,Exact}: 
$$(\pounds_{\vec v} \, g)_{\mu\nu}=2cg_{\mu\nu}
$$
where $\pounds_{\vec v}$ denotes the Lie derivative along $\vec v$.
From standard results \cite{Yano} $(\pounds_{\vec v} \, 
R)^{\alpha}{}_{\beta\gamma\delta}=0$, and as the covariant derivative 
commutes with the Lie derivative along homothetic vectors, then 
$(\pounds_{\vec v} \, \nabla R)_{\mu}{}^{\alpha}{}_{\beta\gamma\delta}=0$. 
Expanding this expression and using on the one hand the 2-symmetry 
and on the other that $\nabla_{\nu}v^{\mu}=c\delta_{\nu}^{\mu}$ one 
derives
$$
(\pounds_{\vec v} \, \nabla R)_{\mu}{}^{\alpha}{}_{\beta\gamma\delta}=
3c\nabla_{\mu}R^{\alpha}{}_{\beta\gamma\delta}=0
$$
from where either $c=0$ or the space is locally symmetric.\finn

\subsection{Conformal and Ricci $k$-symmetry}
\label{conformal+ricci}
All the previous definitions of locally higher-order symmetric spaces
can be straightforwardly generalized by substituting 
the Weyl tensor $C^{\alpha}{}_{\beta\gamma\delta}$ or the Ricci 
tensor $R_{\mu\nu}$ for the Riemann tensor in the defining 
conditions. Thus, one uses the terms {\em locally conformally symmetric} \cite{CG,AM} if 
$$\nabla_{\mu} C^{\alpha}{}_{\beta\gamma\delta}=0,$$
and {\em locally Ricci-symmetric} \cite{Su} if 
$$\nabla_{\rho} R_{\mu\nu}=0.$$ 
Similarly for conformally (or Ricci) recurrent \cite{AM1,Pa,R3,Ol}, 
conformally (or Ricci) 2-recurrent \cite{CR,N,R1}, and 
$k$-recurrent. We can also adopt such a convention and use the terms 
conformal $k$-symmetric and Ricci $k$-symmetric in the obvious way.

Of course, $k$-symmetry implies both conformal $k$-symmetry and 
Ricci $k$-symmetry, but any of the latter does not by itself imply 
the former, see e.g. \cite{R}.
It must be noted that, as above, 
\begin{itemize}
\item proper Riemannian 
Ricci $k$-symmetric spaces are Ricci-symmetric,
\item proper Riemannian 
conformally $k$-symmetric manifolds are conformally symmetric 
\cite{Ta}, and furthermore either locally symmetric or conformally 
flat \cite{DR}.
\end{itemize}

The interplay between the different mentioned 
symmetry and recurrent conditions of any order $k$ has been studied 
in many papers, e.g \cite{Th1,G,Ta,R,R2,DR1,DR2,Hall1,E-K,E-K1}. One 
can also check the Appendix in \cite{DDVV} for an exhaustive list of 
curvature conditions and their overlappings.

\section{Lorentzian 2-symmetry}
\label{lor-bis}
There are several ways to attack the problem of $k$-symmetric and 
$k$-recurrent spaces, among them I would like to cite the following:
\begin{enumerate}
\item pure classical standard tensor 
calculus by using the Ricci and Bianchi identities. This allows us to
obtain necessary restrictions such as those presented already in
subsection \ref{identities};
\item study of 
covariantly constant (also called {\em parallel}) tensor and vector fields, 
and their implications on the manifold local holonomy structure. This
will be dealt with mainly, but not only, in subsection \ref{holonomy};
\item probing the existence of Killing or homothetic vector fields and
consequences thereof (subsection \ref{homothetic} and proof of
proposition \ref{constornull});
\item implications on the curvature invariants and in particular 
their possible constancy or vanishing. This will be largely analyzed
in subsection \ref{invariants}.
\end{enumerate}
Of course, all these methods are clearly interrelated. Here, I 
shall follow a mixed strategy, using results from all of them, trying 
to optimize the path to the sought for results. It turns out that the 
so-called ``superenergy" and/or causal-tensor techniques \cite{S,BS} are 
extremely useful for getting the results as they provide properties
and positive quantities associated with tensors and their tensor
products that can be used 
to, at least partly, replace the ordinary positive-definite metric available in proper 
Riemannian cases.

\subsection{Holonomy and reducibility in Lorentzian manifolds}
\label{holonomy}
To start with, we will need some basic lemmas on local holonomy 
structure. The classical result here is the de Rham decomposition 
theorem \cite{deR,KN} for positive-definite metrics. However, this 
theorem does not hold as such for other signatures, and one has to 
introduce the so-called {\em non-degenerate reducibility} due to Wu 
\cite{Wu,Wu1,Wu2}, who extended de Rham's results to indefinite 
metrics. See also \cite{BI} for the particular case of Lorentzian 
signature. 

To fix the ideas, recall that the holonomy group \cite{KN} 
of $\espaitemps$ is called 
\begin{itemize}
\item {\em reducible} (when acting on the tangent 
spaces) if it leaves a non-trivial subspace of $T_p\varietat$ 
invariant;
\item {\em non-degenerately reducible} if it leaves a 
non-degenerate subspace ---that is, such that the restriction of the 
metric is non-degenerate--- invariant.
\end{itemize}

Fortunately, we will only need a simple result which relates the 
existence of covariantly constant tensors to the holonomy group of 
the manifold in the case of Lorentzian signature. This is a synthesis 
(adapted to our purposes) of the results in \cite{Hall} (see also 
\cite{P,Exact}) but generalized to arbitrary dimension
$n$: \footnote{This result also follows from corollary 4 
and proposition 5 in \cite{SLN}, and can be seen as a consequence of
the results on the holonomy group of irreducible Lorentzian manifolds
found in \cite{SO}.}
\begin{lemma}
\label{redornull}
Let $D\subset \varietat$ be a simply connected domain of an 
$n$-dimensional Lorentzian manifold $\espaitemps$ and assume that 
there exists a non-zero covariantly constant symmetric tensor field 
$h_{\mu\nu}$ not proportional to the metric. Then $(D,g)$ is 
reducible, and further it is not non-degenerately reducible only if 
there exists a null covariantly constant vector field which is the 
unique (up to a constant of proportionality) constant vector field.
\end{lemma}
{\bf Remarks:}
\begin{enumerate}
\item If there is a covariantly constant 1-form $v_{\mu}$, then so is 
obviously $h_{\mu\nu}=v_{\mu}v_{\nu}$ and the manifold (arbitrary 
signature) is reducible, the Span of $v^{\mu}$ being invariant by the 
holonomy group. If $v_{\mu}$ is {\em not} null, then $\espaitemps$ is 
actually {\em non-degenerately} reducible (also called decomposable, 
see \cite{P,Exact}). In this case, the metric can be decomposed into 
two orthogonal parts as $g_{\mu\nu}=c\, v_{\mu}v_{\nu} 
+(g_{\mu\nu}-c\, v_{\mu}v_{\nu})$, where $c=1/(v^{\mu}v_{\mu})$ is 
constant. Thus, necessarily 
$g_{\mu\nu}$ is a {\em flat extension} \cite{RWW} of a 
$(n-1)$-dimensional non-degenerate metric $g_{\mu\nu}-c\, 
v_{\mu}v_{\nu}$.
\item If there is a covariantly constant non-symmetric tensor 
$H_{\mu\nu}$, then its symmetric part is also a constant tensor at 
our disposal, so that one can put $h_{\mu\nu}=H_{(\mu\nu)}$ in the 
lemma unless this vanishes, see next remark.
\item In the case that $H_{\mu\nu}=H_{[\mu\nu]}\neq 0$ is 
antisymmetric, then in fact one can define 
$H_{\mu\rho}H_{\nu}{}^{\rho}=h_{\mu\nu}$, which is symmetric, 
covariantly constant, non-zero and {\em not} proportional to the 
metric if $n>2$. These last two statements follow for example from
lemma 3.2 and corollary 4.1 in \cite{BS}. 
\item Actually, the previous point can be generalized to an arbitrary 
covariantly constant $p$-form $\S_{\mu_1\dots\mu_p}$ (with $p<n$) by defining 
$h_{\mu\nu}=\S_{\mu\rho_2\dots\rho_p}\S_{\nu}{}^{\rho_2\dots\rho_p}$.
\end{enumerate}

\proof The proof relies on the algebraic (or Segre) classification of 2-index 
symmetric covariant tensors, according to its possible eigenvalues and
eigenvectors {\em with respect to the metric}, that is, the solution to the problem $(T_{\mu\nu}-\lambda g_{\mu\nu})b^{\nu}=0$. 
Using the Segre notation, and as is known \cite{O} pp.261-262, only four different types 
occur, namely $[1,1\dots 1], [21\dots 1], [31\dots 1]$ and 
$[z\bar{z}1\dots 1]$, where a comma separates the eigenvalue 
associated to the timelike eigenvector in case this exists, and the 
notation $z\bar{z}$ is used to indicate a pair of complex conjugate 
eigenvalues. Degeneracies are denoted, as usual, by round brackets, 
so for example $[(1,1\dots 1)]$ is the Segre type of the metric 
tensor.

Following \cite{Hall}, one can first of all show that the Segre 
type of $h_{\mu\nu}$ is the same everywhere on $D$ and that its 
eigenvalues are constant on $D$. This follows by computing the 
eigenvalues and eigenvectors at a neighbourhood $U(p)\subset D$, 
which can be taken as differentiable on adequate $U(p)$. Then, for an 
arbitrary curve on $U(p)$ with tangent vector $\vec w$ one can 
parallely transport the eigenvectors to get
$$
h_{\mu\nu}e^{\nu}=\lambda e_{\mu} \hspace{3mm} \Longrightarrow 
\hspace{3mm}
\nabla_{\vec w}(h_{\mu\nu}e^{\nu})=0=(\nabla_{\vec w}\lambda ) e_{\mu}
$$
and given that $\vec w$ is arbitrary one readily gets $\lambda 
=$const.\ on $D$. 

Suppose then that there is an eigenvalue $\lambda$ corresponding to a 
spacelike eigenvector (which is always the case if $n>3$). Suppose 
further that $\lambda$ is not the eigenvalue of the timelike 
eigenvector for case $[1,1\dots 1]$, or of the null eigenvector for 
the cases $[21\dots 1]$ or $ [31\dots 1]$. In other words, by 
relabelling the spacelike eigenvalue if necessary, we are considering 
all cases other than $[(1,1\dots 1)]$, $[(21\dots 1)]$ and $[(31\dots 
1)]$. Then, there is an orthonormal basis of eigenvectors 
$\{\vec{e}_A\}$ for the {\em spacelike} eigenspace corresponding to 
$\lambda$, where $A,B\dots =1,\dots ,q$ and $q$ is the dimension of 
the eigenspace, so that using that $\lambda$ is constant on $D$
$$
h^{\mu}{}_{\nu}e_A^{\nu}=\lambda e_{A}^{\mu} \hspace{3mm} 
\Longrightarrow \hspace{3mm}
h^{\mu}{}_{\nu}\nabla_{\rho} e_A^{\nu}=\lambda \nabla_{\rho} 
e_{A}^{\mu}
\hspace{3mm} \Longrightarrow \hspace{3mm}
\nabla_{\rho} e_A^{\nu}=\sum_B Q_{\rho}{}_{AB} \, e_B^{\nu} 
$$
for some one-forms $Q_{\rho}{}_{AB}$. However, from the 
orthonormal property one easily gets 
$Q_{\rho}{}_{AB}=-Q_{\rho}{}_{BA}$, so that an elementary calculation 
provides
$$
\nabla_{\rho}\left(\sum_A  e_{A}^{\mu}e_{A}^{\nu}\right)=\sum_{AB} 
Q_{\rho}{}_{AB}\left( 
e_{A}^{\mu}e_{B}^{\nu}+e_{B}^{\mu}e_{A}^{\nu}\right)=0
$$
which tells us that the tensor $k_{\mu\nu}=\sum_A  
e_{A}{}_{\mu}e_{A}{}_{\nu}$ is covariantly constant and also, being 
idempotent ($k_{\mu\rho}k^{\rho}{}_{\nu}=k_{\mu\nu}$), a projector. 
It follows that the spacetime is non-degenerately reducible and the 
metric can be non-degenerately decomposed {\em at least} as 
$g_{\mu\nu}=k_{\mu\nu}+(g_{\mu\nu}-k_{\mu\nu})$.

Given that the tensor $h_{\mu\nu}$ is not proportional to the metric, 
the only remaining cases for its Segre type are thus $[(21\dots 1)]$ 
and $[(31\dots 1)]$. Then, a reasoning similar to that in 
\cite{Hall}, as the extra spatial dimensions are immaterial here (but 
the Lorentzian signature is essential), shows that the unique null 
eigendirection can be rescaled to provide a null covariantly constant 
vector field. This must be the only parallel vector field or 
otherwise the previous Remark i would imply that the spacetime is a 
flat extension (ergo non-degenerately reducible) of a 
$(n-1)$-dimensional non-degenerate manifold. (Actually, in the case $[(31\dots 1)]$ there are spacelike parallel vector fields). \finn

\noindent
{\bf Remark:}
Of course, similar reasonings can be applied to the different 
non-degenerate eigenspaces of $h_{\mu\nu}$, providing a more detailed 
descomposition of the spacetime. One should bear in mind this for the rest of the paper together with the basic decomposition $g_{\mu\nu}=k_{\mu\nu}+(g_{\mu\nu}-k_{\mu\nu})$ associated to each non-degenerate eigenspace.

\subsection{Curvature concomitants in 2-symmetric Lorentzian manifolds}
\label{invariants}
Recall that a curvature scalar invariant, see
e.g. \cite{Exact,FKWC,PPCM,P} and references therein, is a 
scalar constructed polynomially from the Riemann tensor, the metric, 
the covariant derivative and possibly the volume element $n$-form of 
$\espaitemps$. They are fundamental in the sense that they are scalars
which depend only
on the metric, and not on the basis used to compute them. They have
been largely studied in General Relativity ($n=4$) as they provide
algebraic classification of spacetimes and important physical
information, see \cite{P,Exact,PPCM} and references therein.
More generally, one can define curvature 1-form
and curvature rank-$r$ tensorial ``invariants'' by using the same ingredients and in the
same manner but leaving 1, \dots , $r$ free indices \cite{FKWC},
respectively. These tensorial quantities are sometimes called 
``curvature concomitants'' (see e.g.\ \cite{Sc} p.15 and p.164), other times ``curvature
covariants'' (e.g.\ \cite{PR} p.260.) Here
we shall use the former.

One can thus  associate three
fundamental natural numbers to all curvature tensorial concomitants, called the
{\em rank}, the {\em degree}, and the {\em order}:
\begin{itemize}
\item the rank is the tensor rank of the tensorial concomitant ---that is,
the number $r$ of free indices---, so that $r=0$ for
scalar invariants.
\item the degree is the (maximum) power of the Riemann tensor used in
 the tensorial concomitant. Hence, they are called linear, quadratic, cubic, etcetera if 
they are linear, quadratic, cubic, and so on, on the Riemann tensor. 
\item finally, the order is the (maximum) number of 
covariant derivatives involved in one single summand of the tensorial concomitant.
\end{itemize}
Curvature tensorial concomitants are said to be homogeneous with respect to the
order if they have the same number of 
covariant derivatives in all its terms, and similarly for homogeneous
with respect to the degree. Of course, all non-homogeneous tensorial concomitants can 
be broken into their respective homogeneous pieces, and therefore in 
what follows we will only consider the homogeneous ones. 

The only linear scalar invariant is $R$ (order zero), and examples of 
quadratic concomitants are $R_{\mu\nu}R^{\mu\nu}$, 
$C_{\alpha\beta\mu\nu}R^{\alpha\beta\mu\nu}$ (order zero), 
$\nabla_{\mu}R_{\alpha\beta\gamma\delta}\nabla^{\mu}R^{\alpha\beta\gamma\delta}$ 
(order 2), 
$R_{\alpha\beta\mu\nu}R^{\alpha\beta\mu\nu}+\nabla_{\rho}R_{\mu\nu}\nabla^{\rho}R^{\mu\nu}$ 
(non-homogeneous), $R\nabla_{\mu} R$ (order 1, rank 1), 
$\nabla^{\rho}C^{\alpha\beta\mu\nu}\nabla_{\rho}R_{\beta\nu}$ (order 
2, rank 2). 

The following are very useful lemmas.
\begin{lemma}
\label{nullorzero}
Let $D\subset \varietat$ be a simply connected domain and let $g$ have arbitrary signature. Any 
curvature 1-form concomitant which is covariantly constant must be 
necessarily null (possibly zero).
\end{lemma}
\proof If it were non-null, due to the Remark \rm{i} of Lemma 
\ref{redornull}, $(D,g)$ would be a flat extension of a $(n-1)$ 
space, but this is clearly impossible as the curvature would then be 
that of only the reduced $(n-1)$ part.\finn
\begin{lemma}
\label{2form}
Let $(D,g)$ be as in lemma \ref{redornull}. Any
curvature $2$-form ($2<n$) concomitant which is covariantly constant is null
(possibly zero), and if non-zero there is a null covariantly
constant vector field.
\end{lemma}
(For the definition of null $2$-form, see e.g.\ \cite{Exact,BS})

\proof From the Remark \rm{iii} of Lemma 
\ref{redornull}, and using the notation there, $(D,g)$ is reducible. 
If the $2$-form is null and non-zero, then
in fact the parallel tensor $h_{\mu\nu}=k_{\mu}k_{\nu}$ where $k_{\mu}$
defines the unique null direction determined by the null
$2$-form. Therefore in this case $k_{\mu}$ must be covariantly
constant. If the $2$-form were non-null, then $(D,g)$ would be
non-degenerately reducible so that the 2-form concomitant itself
decomposes as the sum of the respective reduced 2-forms, some of which
may be zero. The non-zero ones would necessarily be such that
the corresponding reduced $h_{\mu\nu}$ is proportional to the 
metric of a 2-dimensional subspace defined by the decomposition. As 
this 2-space is necessarily 2-symmetric, by standard results it must 
in fact be of constant curvature and thus the curvature 2-form concomitant
would have to vanish anyway. \finn

A very simple yet useful lemma is the following.
\begin{lemma}
\label{orderdegree}
Let  $\espaitemps$ be an $n$-dimensional {\em 2-symmetric} Lorentzian
manifold. Then
\begin{itemize}
\item all non-zero curvature tensorial concomitants have a degree 
necessarily greater or equal than the order;
\item all curvature tensorial concomitants with order equal to degree are
necessarily parallel.
\end{itemize}
\end{lemma}
\proof As the spacetime is 2-symmetric, there will be no non-zero
curvature tensorial concomitants 
involving derivatives of the Riemann tensor of order higher than one. This proves the first
assertion and, mutatis mutandis, the second too.\finn

With these lemmas in hand, one of the main results can be proven. 
Observe that from lemma \ref{nullorzero} it follows that, in 2-symmetric 
spaces, either $R$ is constant or $\nabla_{\mu} R$ is null and 
covariantly constant. This is a particular example of the following 
general result.
\begin{prop}
\label{constornull}
Let $D\subset \varietat$ be a simply connected domain of an 
$n$-dimensional {\em 2-symmetric} Lorentzian manifold $\espaitemps$. 
Then either
\begin{itemize}
\item all scalar invariants of the Riemann tensor of 
order $m$ and degree up to $m+2$ are constant on $D$; or
\item there is a covariantly constant null vector field on $D$.
\end{itemize}
\end{prop}
\proof Let $I$ be a scalar invariant of order $m$. If its degree is also 
$m$ then obviosuly (see Lemma \ref{orderdegree})
$\nabla I =0$ so that $I$ is constant. If its 
degree is $m+1$, then $\nabla_{\nu}\nabla_{\mu} I =0$ from where 
$v_{\mu}=\nabla_{\mu} I$ is a covariantly constant curvature 1-form
concomitant, so that from Lemma \ref{nullorzero} it must be null or 
zero, in the last case implying that $I$ is constant. Finally, if the 
degree is $m+2$, then $\nabla_{\rho}\nabla_{\nu}\nabla_{\mu} I =0$ so 
that $h_{\mu\nu}=\nabla_{\nu}\nabla_{\mu} I=\nabla_{\nu}v_{\mu}$ is a 
covariantly constant symmetric tensor. According to Lemma 
\ref{redornull}, there arise three possibilities: there is a null 
covariantly constant vector field, or $(D,g)$ is non-degenerately 
reducible and there is no null parallel vector field, or 
$h_{\mu\nu}=cg_{\mu\nu}$ for a constant $c$. In the first possibility 
we are done; in the second case, the space is decomposed into several 
irreducible and necessarily 2-symmetric parts, in particular
one for which 
$h_{\mu\nu}$ is proportional to its reduced metric tensor. 
As $I$ will be a function 
of the corresponding scalar invariants of the irreducible parts,
everything reduces effectively to only the third possibility. But then
Lemma \ref{homo} implies that either $c=0$ or the (reduced) space is
locally symmetric. In both cases $h_{\mu\nu}=\nabla_{\nu}v_{\mu}=0$ implying that 
$v_{\mu}=\nabla_{\mu} I$ must be covariantly constant. Thus, the 
same reasoning as before using Lemma \ref{nullorzero} ends the proof. 
\finn

There are, of course, ways to obtain particular results of this lemma 
in a direct manner. For example, by contracting $\alpha$ and 
$\gamma$, $\beta$ and $\delta$, $\nu$ and $\mu$ in (\ref{RicDR}) one 
gets $R^{\rho}{}_{\mu}\nabla_{\rho}R=0$, so that contracting here 
with $\nabla^{\mu}$ and using the contracted Bianchi identity
\be
\nabla_{\mu}R^{\mu}{}_{\nu}=\frac{1}{2}\nabla_{\nu}R 
\label{contBianchi2}
\ee
the condition $\nabla_{\nu}R\nabla^{\nu}R=0$ is obtained. Then either 
$R=$const.\ or $\nabla_{\nu}R$ is covariantly constant and null. 
Nevertheless, there are other results which are much more difficult 
to prove using pure tensor calculus. For instance, 
\be
\nabla_{\alpha}(R^{\mu\nu}R_{\mu\nu})=0= 
R^{\mu\nu}\nabla_{\alpha}R_{\mu\nu} \label{R^2} 
\ee
which is obviously contained in Proposition \ref{constornull} if 
there is no null parallel vector field. 

The previous proposition has immediate consequences providing more 
information about curvature tensorial concomitants. For instance
\begin{coro}
\label{consequences}
Under the conditions of Proposition \ref{constornull}, either there 
is a covariantly constant null vector field on $D$ or the following 
statements hold
\begin{enumerate}
\item All curvature scalar invariants of any order and degree formed 
as functions of the homogeneous ones of order $m$ and degree up to 
$m+2$ are constant on $D$;
\item All curvature 1-form concomitants of order $m$ and degree up to 
$m+1$ are zero.
\item All curvature scalar invariants with order equal to degree vanish.
\item All curvature rank-2 tensorial concomitants with order equal to degree are 
zero.
\end{enumerate}
\end{coro}
{\bf Remark:} Of course, it can happen that the mentioned curvature 
concomitants vanish {\em and} there is a covariantly constant null 
vector field too.

\noindent
\proof The first statement is trivial. 
For the second, choose any 1-form concomitant $I_{\mu}$ of order $m$. 
If the degree is also $m$, then obviously (Lemma \ref{orderdegree})
$I_{\mu}$ is covariantly 
constant so that Lemma \ref{nullorzero} applies. If its degree is 
$m+1$,
then its exterior differential $F_{\mu\nu}=\nabla_{[\nu}I_{\mu]}$ is 
a covariantly constant 2-form and thus Lemma 
\ref{2form} implies that, if there is no null parallel vector field (which 
means that $F_{\mu\nu}$ cannot be non-zero and null), 
$\nabla_{[\nu}I_{\mu]}= 0$, so
that $I_{\mu}$ is an exact 1-form and therefore $I_{\mu}$ is the
differential of a function $f$, $I_{\mu}=\nabla_{\mu}f$. If $f$ happens to be a scalar invariant $f=I$, then $I$ is of order $m-1$ and degree $m+1$. 
As $I$ is constant due to Proposition \ref{constornull}, then 
$I_{\mu}=0$. But even in the case that $f$ is not a scalar curvature invariant, the symmetric tensor $h_{\mu\nu}=\nabla_{\nu}I_{\mu}=\nabla_{\nu}\nabla_{\mu}f$ is covariantly constant, and a reasoning totally analogous to that in the end of the proof of Proposition \ref{constornull} implies that $I_{\mu}$ is parallel, so that Lemma \ref{nullorzero} applies again.

For the third statement, take any scalar invariant $I$ with order and 
degree equal to $m$. Then obviously this is the divergence 
$I=\nabla_{\mu}I^{\mu}$ of a 1-form concomitant $I_{\mu}$ of order 
$m-1$ and degree $m$. As $I_{\mu}=0$ due to point {\it ii}, $I=0$.

Finally, for the point $iv$, as the curvature rank-2 tensorial concomitant 
$I_{\mu\nu}$ has the order equal to the degree, it is covariantly 
constant $\nabla_{\rho}I_{\mu\nu}=0$. Furthermore, due to point {\it 
iii}, its trace vanishes $I^{\nu}{}_{\nu}=0$. Its symmetric part 
$I_{(\mu\nu)}$ is therefore covariantly constant and {\em not} 
proportional to the metric. By using Lemma \ref{redornull}, following 
the same reasoning as in Proposition \ref{constornull} one easily 
proves $I_{(\mu\nu)}=0$ unless there is a null parallel vector field. 
If this is not the case, then $I_{\mu\nu}=I_{[\mu\nu]}$ is a 
covariantly constant 2-form and then lemma \ref{2form} implies 
finally that $I_{\mu\nu}=0$.\finn

There is a very long list of vanishing curvature scalar and tensorial concomitants as a 
result of this proposition---if there is no null parallel vector 
field---. The following is the list of the quadratic ones (only an 
independent set \cite{FKWC}, apart from (\ref{R^2}), is given, also 
omitting those contaning $\nabla_{\mu}R=0$):
\bea
R^{\mu\nu}\nabla_{\mu}R_{\nu\alpha}=0, 
R^{\mu\nu\rho\alpha}\nabla_{\mu}R_{\nu\rho}=0,\,\,
R^{\mu\nu\rho\sigma}\nabla_{\mu}R_{\nu\rho\sigma\alpha}=0=
R^{\mu\nu\rho\sigma}\nabla_{\alpha}R_{\mu\nu\rho\sigma}, \label{r1}\\
\nabla_{\alpha}R^{\mu\nu}\nabla_{\beta}R_{\mu\nu}=
\nabla_{\mu}R_{\nu\beta}\nabla_{\alpha}R^{\mu\nu}=
\nabla_{\mu}R_{\nu\alpha}\nabla^{\mu}R^{\nu}{}_{\beta}=
\nabla_{\mu}R_{\nu\alpha}\nabla^{\nu}R^{\mu}{}_{\beta}=0, 
\label{r2ric}\\
\nabla^{\mu}R^{\nu\rho}\nabla_{\alpha}R_{\beta\rho\mu\nu}=
\nabla^{\mu}R^{\nu\rho}\nabla_{\mu}R_{\alpha\nu\beta\rho}=0,\,\,
\nabla_{\alpha}R^{\mu\nu\rho\sigma}\nabla_{\beta}R_{\mu\nu\rho\sigma}=
\nabla^{\sigma}R^{\mu\nu\rho\alpha}\nabla_{\sigma}R_{\mu\nu\rho\beta}=0\label{r2rie}
\eea
where of course the traces of (\ref{r2ric}-\ref{r2rie}) vanish, and 
one can also write the same expressions using the Weyl tensor 
instead of the Riemann tensor.

\section{Main results}
\label{main}
All necessary results to prove the main theorems have now been 
gathered. Let us start with an important result which will be derived 
by using the so-called causal tensors and ``superenergy" techniques 
\cite{S,BS}. The idea is that a tensor such as 
$\nabla_{\mu}R_{\alpha\beta}$, whose tensor square has all double 
traces vanishing---see (\ref{r2ric})---, and from (\ref{treq4}) also 
satisfies 
\be
(\nabla_{\nu}R^{\rho}{}_{\mu}-\nabla_{\mu}R^{\rho}{}_{\nu})
\nabla_{\rho}R_{\alpha\beta}=0, \,\,\, 
(\nabla^{\rho}R_{\mu\nu}-2\nabla_{\nu}R^{\rho}{}_{\mu})
\nabla_{\rho}R_{\alpha\beta}=0, \label{treq4Ric}
\ee
can only be non-zero if $\nabla_{\mu}R_{\alpha\beta}=
k_{\mu}(k_{\alpha}p_{\beta}+k_{\beta}p_{\alpha})$ where $k_{\mu}$ is
a null vector field orthogonal to $p_{\mu}$. But as the tensor 
$\nabla_{\mu}R_{\alpha\beta}$ is covariantly constant, so is the null 
vector field.
\begin{theorem}
\label{ric-flat}
Let $D\subset \varietat$ be a simply connected domain of an 
$n$-dimensional {\em 2-symmetric} Lorentzian manifold $\espaitemps$. 
Then, if there is no null covariantly constant vector field on $D$, 
$(D,g)$ is the direct product of irreducible submanifolds each of which is either Ricci-flat ($R_{\mu\nu}=0$) or locally symmetric.
\end{theorem}
\proof The first part of the reasoning is algebraic and can be 
performed at any point $p\in D$. Take the double (2,1)-form 
$\nabla_{[\alpha}R_{\beta]\lambda}$ and construct its basic 
superenergy tensor \cite{S}:
\bean
T_{\alpha\beta\lambda\mu}=(\nabla_{\alpha}R_{\rho\lambda}-\nabla_{\rho}R_{\alpha\lambda})
(\nabla_{\beta}R^{\rho}{}_{\mu}-\nabla^{\rho}R_{\beta\mu})+
(\nabla_{\alpha}R_{\rho\mu}-\nabla_{\rho}R_{\alpha\mu})
(\nabla_{\beta}R^{\rho}{}_{\lambda}-\nabla^{\rho}R_{\beta\lambda})-\\
-\frac{1}{2} 
g_{\alpha\beta}(\nabla_{\sigma}R_{\rho\lambda}-\nabla_{\rho}R_{\sigma\lambda})
(\nabla^{\sigma}R^{\rho}{}_{\mu}-\nabla^{\rho}R^{\sigma\mu})-
g_{\lambda\mu}(\nabla_{\alpha}R_{\rho\sigma}-\nabla_{\rho}R_{\alpha\sigma})
(\nabla_{\beta}R^{\rho\sigma}-\nabla^{\rho}R_{\beta}{}^{\sigma})+\\
+\frac{1}{4} g_{\alpha\beta} g_{\lambda\mu}
(\nabla_{\tau}R_{\rho\sigma}-\nabla_{\rho}R_{\tau\sigma})
(\nabla^{\tau}R^{\rho\sigma}-\nabla^{\rho}R^{\tau\sigma}) . 
\hspace{2cm}
\eean
This is the (essentially unique) tensor which is quadratic on 
$\nabla_{[\mu}R_{\alpha]\beta}$ and is ``causal" or ``future" 
\cite{S,BS}, that is to say, it satisfies the {\em dominant property}:
$$
T_{\alpha\beta\lambda\mu}u_1^{\alpha}u_2^{\beta}u_3^{\lambda}u_4^{\mu}\geq 
0
$$
for every choice of future-pointing vectors 
$\{u_1^{\alpha},u_2^{\beta},u_3^{\lambda},u_4^{\mu}\}$ (as a matter 
of fact, this is strictly positive if all the previous vectors are 
timelike). It is known that the tensor $T_{\alpha\beta\lambda\mu}$ 
vanishes if and only if so does $\nabla_{[\mu}R_{\alpha]\beta}$, and 
if and only if its contraction in all indices with a timelike vector  
vanishes. Another  direct property is 
$T_{\alpha\beta\lambda\mu}=T_{(\alpha\beta)(\lambda\mu)}$. In our 
case, first of all by use of (\ref{r2ric}) the tensor becomes 
$$
T_{\alpha\beta\lambda\mu}=(\nabla_{\alpha}R_{\rho\lambda}-\nabla_{\rho}R_{\alpha\lambda})
(\nabla_{\beta}R^{\rho}{}_{\mu}-\nabla^{\rho}R_{\beta\mu})+
(\nabla_{\alpha}R_{\rho\mu}-\nabla_{\rho}R_{\alpha\mu})
(\nabla_{\beta}R^{\rho}{}_{\lambda}-\nabla^{\rho}R_{\beta\lambda})
$$
and using repeatedly (\ref{treq4Ric}) this finally simplifies to
\be
T_{\alpha\beta\lambda\mu}=\nabla_{\alpha}R_{\rho\lambda}\nabla_{\beta}R^{\rho}{}_{\mu}+
\nabla_{\alpha}R_{\rho\mu}\nabla_{\beta}R^{\rho}{}_{\lambda} . 
\label{se}
\ee
Now, one has to use several elementary properties of causal tensors. 
As $T_{\alpha\beta\lambda\mu}$ is future (theorem 4.1 in \cite{S}), 
its contraction in any number of indices with arbitrary 
future-pointing vectors is also a causal tensor (property 2.2 in 
\cite{BS}), so that in particular, for an arbitrary
timelike vector $u^{\mu}$, one has that
\be
T_{\alpha\beta\lambda\mu}u^{\alpha}u^{\beta}=
2\, u^{\sigma}\nabla_{\sigma}R_{\rho\lambda}\, 
u^{\tau}\nabla_{\tau}R^{\rho}{}_{\mu} \label{lon1}
\ee
is also a future tensor, non-vanishing unless
$\nabla_{[\mu}R_{\nu]\rho}=0$. 
But this means\footnote{Observe the change of
signature with respect to Ref.\cite{BS}, which may be a little
confusing. One should not use, for instance, proposition 2.1 in
\cite{BS}, but rather the corollary 2.1 in that paper.}, (corollaries
2.1 or 2.7 in \cite{BS}) that {\em for arbitrary timelike $u^{\rho}$} the tensor 
\be
\stackrel{u}{T}_{\lambda\mu}\equiv u^{\sigma}\nabla_{\sigma}R_{\lambda\mu}
\label{tu}
\ee
cannot be causal, neither future nor past, unless it takes the form
$k_{\mu}k_{\nu}$ for some null $k^{\mu}$ possibly depending on
$\vec u$. There arise two possibilities.

\vspace{3mm}
\noindent
(i) If $\stackrel{u}{T}_{\lambda\mu} = k_{\lambda}k_{\mu}$ for some $\vec u$,
then (\ref{lon1}) implies that
$T_{\alpha\beta\lambda\mu}u^{\alpha}u^{\beta}=0$, and this leads
(corollary 2.3 in \cite{BS}) to $T_{\alpha\beta\lambda\mu}=0$ which is equivalent 
(property 3.4 in \cite{S}) to $\nabla_{[\nu}R_{\lambda]\mu}=0$. Thus, 
the tensor $\nabla_{\nu}R_{\lambda\mu}$ is completely
symmetric. Contracting (\ref{se}) with $u^{\alpha}v^{\beta}$, where
$v^{\beta}$ is any other timelike vector, we derive
$$
k_{\rho}\, \stackrel{v}{T}\!{}^{\rho}{}_{(\lambda}k_{\mu)}=0 \,\,\,
\Longrightarrow \,\,\, k_{\rho}\, \stackrel{v}{T}\!{}^{\rho}{}_{\lambda}=0
$$ 
for any such $\vec v$.
However, from (\ref{r2ric}) one deduces that, for arbitrary timelike $\vec w$ and $\vec v$
\be
\stackrel{w}{T}_{\lambda\mu} \stackrel{v}{T}\!{}^{\lambda\mu}=
 w^{\sigma}\nabla_{\sigma}R_{\lambda\mu}v^{\tau}\nabla_{\tau}R^{\lambda\mu}=0
\label{uv}
\ee 
which together with the previous implies that, in fact, there is a fixed null 
1-form $k_{\mu}$ independent of $u^{\rho}$ such that, using the complete symmetry of the tensor
\be
\nabla_{\nu}R_{\lambda\mu}=f k_{\nu}k_{\lambda}k_{\mu} .
\label{kkk}
\ee
This reasoning has been performed at a fixed point. One cannot
therefore extract conclusions by taking derivatives of the previous expression
until the second possibility has been analyzed.

\vspace{3mm}
\noindent 
(ii) Suppose then that, at the given point, all the tensors (\ref{tu})
are non-causal, for arbitrary timelike $\vec u$. Furthermore, from
$$
T_{\alpha\beta\lambda\mu}u^{\alpha}u^{\beta}v^{\lambda}v^{\mu}\geq 0
$$
it follows that $\stackrel{u}{T}_{\rho\lambda}v^{\lambda}$ is spacelike for
arbitrary timelike future directed $\vec u,\vec v$. This immediately
implies that  $\stackrel{u}{T}_{\rho\lambda}$ cannot be of algebraic type
$[1,1\dots 1]$ for any $\vec u$ as this type has a timelike
eigenvector. Similarly, $\stackrel{u}{T}_{\rho\lambda}$ cannot be of type
 $[z\bar z 1\dots 1]$ 
for any $\vec u$, because if it were one could easily check that the only way 
$\stackrel{u}{T}_{\rho\lambda}v^{\lambda}$ is spacelike  for
arbitrary timelike future directed $\vec v$ is that
$$
\stackrel{u}{T}_{\mu\nu}=\nu \left(k_{\mu}k_{\nu}-\ell_{\mu}\ell_{\nu}\right)
+\stackrel{u}{P}_{\mu\nu} , \,\,\,\,\, \nu(u)\neq 0
$$
for null $k_{\mu}$ and $\ell_{\nu}$---possibly depending on $\vec
u$--- such that $k^{\mu}\ell_{\mu}=-1$ and with $\stackrel{u}{P}_{\mu\nu}=\stackrel{u}{P}_{(\mu\nu)}$,
$\stackrel{u}{P}_{\mu\nu}\ell^{\mu}=0$,
$\stackrel{u}{P}_{\mu\nu}k^{\mu}=0$. But then $\stackrel{u}{T}_{\mu\nu}\ell^{\mu}$ and
$\stackrel{u}{T}_{\mu\nu}k^{\mu}$ would be both null and 
$$
T_{\alpha\beta\lambda\mu}u^{\alpha}u^{\beta}\ell^{\lambda}\ell^{\mu}=0,\,\,
T_{\alpha\beta\lambda\mu}u^{\alpha}u^{\beta}k^{\lambda}k^{\mu}=0
$$
which would imply (property 2.3 in \cite{BS}) that in fact
$$
T_{\alpha\beta\lambda\mu}\ell^{\lambda}\ell^{\mu}=0, \,\,
T_{\alpha\beta\lambda\mu}k^{\lambda}k^{\mu}=0.
$$
From here one would easily deduce that, {\em for all} timelike $\vec
u$, $\stackrel{u}{T}_{\mu\nu}\ell^{\mu}$ would be null and proportional to $k_{\nu}$
and $\stackrel{u}{T}_{\mu\nu}k^{\mu}$ would be null and proportional to
$\ell_{\nu}$, with both $k_{\nu}$ and $\ell_{\nu}$ independent of $u^{\mu}$. This would lead to 
$$
\nabla_{\rho}R_{\mu\nu}=w_{\rho}\left(k_{\mu}k_{\nu}-\ell_{\mu}\ell_{\nu}\right)+
P_{\rho\mu\nu}
$$
for fixed $w_{\rho}$ and null $k_{\mu}$ and $\ell_{\mu}$ 
with $P_{\rho\mu\nu}=P_{\rho(\mu\nu)}$ and $P_{\rho\mu\nu}\ell^{\mu}=0$,
$P_{\rho\mu\nu}k^{\mu}=0$. Using here that
$\nabla_{\rho}R^{\rho}{}_{\nu}=0$ and (\ref{r2ric}) this would
inevitably lead to $\nabla_{\nu}R_{\lambda\mu}=0$.

The only remaining possibility is that all tensors (\ref{tu}) have
algebraic type $[21\dots 1]$ or $[31\dots 1]$, in both cases having a
null eigenvector $k_{\mu}$---possibly depending on $\vec u$. The
former type $[21\dots 1]$ is easily ruled out by using again the
property (\ref{r2ric}) which leads to $\stackrel{u}{T}_{\lambda\mu} =
k_{\lambda}k_{\mu}$
for null $k_{\mu}$, that is, to the already considered
possibility (i). Let us finally consider the remaining type 
$[31\dots  1]$. With the help once again of (\ref{r2ric}) it follows
that $\stackrel{u}{T}_{\mu\nu}=k_{(\mu}\stackrel{u}{q}_{\nu)}$ with 
$\stackrel{u}{q}_{\mu}$ spacelike and
orthogonal to $k_{\mu}$. But then, 
$T_{\alpha\beta\lambda\mu}u^{\alpha}u^{\beta}k^{\lambda}=0$ so that
$T_{\alpha\beta\lambda\mu}k^{\lambda}=0$ and using the same
type of reasonings as before
$\stackrel{u}{T}_{\mu\nu}=k_{(\mu}\stackrel{u}{q}_{\nu)}$
with a fixed null $k_{\mu}$ which does not depend on $u^{\rho}$. This implies that 
$$
\nabla_{\nu}R_{\lambda\mu}=P_{\nu(\lambda}k_{\mu)}
$$
for some $P_{\nu\lambda}$ with $k^{\lambda}P_{\nu\lambda}=0$. Using now once more the several
formulas in (\ref{r2ric}) one derives the result that, in fact, 
$$
\nabla_{\nu}R_{\lambda\mu}=k_{\nu}q_{(\lambda}k_{\mu)}
$$  
with $q^{\mu}$ spacelike and orthogonal to $k_{\mu}$.

Summarizing the two possibilities (i) and (ii), at every point of the
domain $D$ the
covariant derivative of the Riemann tensor must take the form
$$
\nabla_{\nu}R_{\lambda\mu}=k_{\nu}p_{(\lambda}k_{\mu)}
$$
where $k_{\mu}$ is a differentiable null vector field and $p^{\mu}$ is
a differentiable vector field
orthogonal to $k_{\mu}$ ---it may be proportional to $k_{\mu}$ in some
places, or even zero---. But then
$$
T_{\alpha\beta\lambda\mu}=(p^{\rho}p_{\rho})k_{\alpha}k_{\beta}k_{\lambda}k_{\mu}
$$
so that, as $T_{\alpha\beta\lambda\mu}$ is a parallel tensor from its
definition, it follows that there is a vector proportional to $k^{\mu}$
which is parallel, unless $p_{\mu}$ is null or zero. But in these
cases $T_{\alpha\beta\lambda\mu}=0$ implying
$\nabla_{[\nu}R_{\lambda]\mu}=0$ and therefore (\ref{kkk}) holds. But
then again, as $\nabla_{\alpha}\nabla_{\nu}R_{\lambda\mu}=0$, it easily follows 
that, if $f$ is non-zero, then a null vector field proportional to
$k_{\mu}$ is covariantly constant. 

To finalize, as the existence of a null covariantly constant vector
field is against the hypothesis of the theorem, it necessarily follows
$$
\nabla_{\nu}R_{\lambda\mu}=0 .\label{constRic}
$$
Therefore, the Ricci tensor is a parallel symmetric tensor field, so that 
from Lemma \ref{redornull} the manifold is either non-degenerately 
decomposable or the Ricci tensor is proportional to the metric. By 
decomposing into irreducibe parts, all of them 2-symmetric as has been already explained several times, the only relevant cases are that of Einstein spaces (for the reduced spaces, if necessary)
$$
R_{\lambda\mu}=\frac{R}{m}g_{\lambda\mu} .
$$
(Here $m\leq n$ is the dimension of the given irreducible submanifold if this is necessary, and then $g_{\lambda\mu}$, $R_{\lambda\mu}$ and $R$ stand for its corresponding reduced metric tensor, Ricci tensor and scalar curvature, respectively.)
However, from  this relation and the contracted Bianchi identity ---applied to the reduced parts--- one deduces
$$
\nabla_{\rho}R^{\rho}{}_{\beta\gamma\delta}=0 ,
$$
hence, contracting $\nu$ with $\lambda$ in the equation (\ref{eq4}) 
applied to the corresponding Riemann tensor, it follows that
$$
R\, \nabla_{\mu}R_{\alpha\beta\gamma\delta}=0
$$
so that finally each of the irreducible components of the spacetime is either locally symmetric or its Ricci tensor vanishes, which finishes the proof.\finn

{\bf Remark}. It must be stressed that this proof is only valid for 
Lorentzian manifolds, as the definition of future tensors requires 
this signature. I would like to stress, in passing, that this proof 
shows the mathematical potentialities of the superenergy 
construction and of the theory of causal tensors. 

Now, one can at last prove that the narrow space left between locally 
symmetric and 2-symmetric Lorentzian manifolds can only be filled by 
spaces with a covariantly constant null vector field. In order to
alleviate the path to the final result, we first present the algebraic
part of the reasoning as separate lemmas, interesting on their own.
\begin{lemma}
\label{alg}
At any point of a Ricci-flat Lorentizan manifold $\espaitemps$ where
the Riemann (or equivalently the Weyl)
tensor satisfies the first in (\ref{basicR})
\be
\nabla_{(\tau}C^{\rho}{}_{\nu)\lambda\mu}\nabla_{\rho}C_{\alpha\beta\gamma\delta}=0
\label{basicC}
\ee
there exists a null vector $k^{\mu}$ such that
$$
k^{\rho}\nabla_{\rho}C_{\alpha\beta\gamma\delta}=0,\,\,\,\,
k^{\alpha}\nabla_{\rho}C_{\alpha\beta\gamma\delta}=0.
$$
\end{lemma}
\proof First of all, there must be at least a
tangent vector $p^{\mu}$ such that 
$$
p^{\rho}\, \nabla_{\rho}C_{\alpha\beta\gamma\delta}=0
$$
as otherwise it would follow from (\ref{basicC}) that
$\nabla_{(\tau}C^{\rho}{}_{\nu)\lambda\mu}=0$ which together with the identity
$\nabla_{[\tau}C_{\rho\nu]\lambda\mu}=0$ would lead to
$\nabla_{\tau}C_{\rho\nu\lambda\mu}=0$
anyway. Let us define, for any arbitrary 2-form $F_{\mu\nu}$, the
following tensor
$$
\stackrel{F}{\Omega}\!{}_{\rho\alpha\beta}\equiv 
\nabla_{\rho}C_{\alpha\beta\gamma\delta}F^{\gamma\delta}
$$
which has the immediate properties
$$
\stackrel{F}{\Omega}\!{}_{\rho\alpha\beta}=\stackrel{F}{\Omega}\!{}_{\rho[\alpha\beta]},\,\,\,\,
\stackrel{F}{\Omega}\!{}_{[\rho\alpha\beta]}=0\,\,\,\,\, 
\stackrel{F}{\Omega}\!{}^{\rho}{}_{\rho\beta}=0\, .
$$
The main condition (\ref{basicC}) can then be reexpressed as, for {\em
  any} 2-forms $F^{\mu\nu}$ and $G^{\mu\nu}$,
\be
\stackrel{F}{\Omega}\!{}_{(\tau}{}^{\rho}{}_{\nu)}\stackrel{G}{\Omega}\!{}_{\rho\alpha\beta}=0.
\label{basicC2}
\ee
If the vector $p^{\mu}$ is non-null, it follows directly by orthogonal decomposition
\be
\stackrel{F}{\Omega}\!{}_{\rho\alpha\beta}=\stackrel{F}{T}\!{}_{\rho\alpha\beta}+2\,
\stackrel{F}{C}\!{}_{\rho[\beta}p_{\alpha]} \label{OT}
\ee
where
\bean
\stackrel{F}{T}\!{}_{\rho\alpha\beta}=\stackrel{F}{T}\!{}_{\rho[\alpha\beta]},\,\,\,\, 
\stackrel{F}{T}\!{}_{[\rho\alpha\beta]}=0,\,\,\,\,
\stackrel{F}{T}\!{}^{\rho}{}_{\rho\beta}=0,\,\,\,\,
\stackrel{F}{C}\!{}_{\rho\beta}=\stackrel{F}{C}\!{}_{(\rho\beta)},\,\,\,\,
\stackrel{F}{C}\!{}^{\rho}{}_{\rho}=0,\\
\stackrel{F}{T}\!{}_{\rho\alpha\beta}p^{\rho}=0, \,\,\,\,\,
\stackrel{F}{T}\!{}_{\rho\alpha\beta}p^{\alpha}=0,\,\,\,\,\,
\stackrel{F}{C}\!{}_{\rho\beta}p^{\rho}=0.\hspace{2cm}
\eean
But then, it is very easy to see that (\ref{basicC2}) implies the
following conditions for arbitrary $F^{\mu\nu}$ and $G^{\mu\nu}$
\bean
\stackrel{F}{C}\!{}_{\tau}{}^{\rho}\stackrel{G}{C}\!{}_{\rho\alpha}&=&0,\\
\stackrel{F}{T}\!{}_{(\tau}{}^{\rho}{}_{\nu)}\stackrel{G}{C}\!{}_{\rho\alpha}=0,\,\,\,\,\,\,\,\,\,
\stackrel{F}{C}\!{}_{\tau}{}^{\rho} \stackrel{G}{T}\!{}_{\rho\alpha\beta}=0,\,\,\,\,\,
&\Longrightarrow& \,\,\,\, \stackrel{F}{C}\!{}_{\tau}{}^{\rho}
\stackrel{G}{T}\!{}_{\alpha\rho\beta}=0\\
\stackrel{F}{T}\!{}_{(\tau}{}^{\rho}{}_{\nu)} \stackrel{G}{T}\!{}_{\rho\alpha\beta}&=&0
\eean
where the implication in the middle line follows on using the property 
$\stackrel{F}{T}\!{}_{[\rho\alpha\beta]}=0$. The first of these conditions
implies necessarily (by repeated application of corollary 2.7 in \cite{BS}) that
$$
\stackrel{F}{C}\!{}_{\rho\alpha}=f(F)k_{\rho}k_{\alpha}
$$
for a {\em fixed} null vector $k_{\mu}$ independent of $F\in \Lambda_2$
and a function---which may
have zeros--- depending on $F\in \Lambda_2$. But then the second
condition implies that, for any $F,G\in \Lambda_2$,
$$
f(F)\, \stackrel{G}{T}\!{}_{\rho\alpha\beta}k^{\rho}=0,\,\,\,
f(F)\, \stackrel{G}{T}\!{}_{\rho\alpha\beta}k^{\alpha}=0
$$
so that either $f(F)=0$ or all the tensors
$\stackrel{G}{T}\!{}_{\rho\alpha\beta}$ are completely orthogonal to $k^{\rho}$
for arbitrary $G\in \Lambda_2$. In the latter case, it follows that 
$k^{\rho}\nabla_{\rho}C_{\alpha\beta\gamma\delta}=0$. In the
former case, all tensors $\stackrel{F}{C}\!{}_{\rho\alpha}$ vanish and therefore
one can deduce, by using the symmetry and trace properties of 
$\nabla_{\rho}C_{\alpha\beta\gamma\delta}$, that this tensor is
orthogonal in all its indices to $p^{\mu}$. But then one can start
again the reasoning now restricted to the tangent subspace 
$\{\vec p\}^{\perp}$ by choosing another $\vec{p}\,'\in \{\vec
p\}^{\perp}$ such that
$p'^{\rho}\nabla_{\rho}C_{\alpha\beta\gamma\delta}=0$, and so on until
either there is a null $k^{\mu}$ such that 
$k^{\rho}\nabla_{\rho}C_{\alpha\beta\gamma\delta}=0$ or the tensor
$\nabla_{\rho}C_{\alpha\beta\gamma\delta}$ vanishes. 

All in all, the conclusion is that we could have
assumed, from the beginning, that $p^{\rho}=k^{\rho}$ is null. In this
case, one can directly write, analogously to (\ref{OT}),
$$
\stackrel{F}{\Omega}\!{}_{\rho\alpha\beta}=\stackrel{F}{S}\!{}_{\rho\alpha\beta}+2\,
\stackrel{F}{A}\!{}_{\rho[\beta}\ell_{\alpha]}
$$
where $\ell_{\mu}$ is null such that $k^{\mu}\ell_{\mu}=-1$ and 
\bean
\stackrel{F}{S}\!{}_{\rho\alpha\beta}=\stackrel{F}{S}\!{}_{\rho[\alpha\beta]},\,\, 
\stackrel{F}{S}\!{}_{[\rho\alpha\beta]}=0,\,\,
\stackrel{F}{A}\!{}_{\rho\beta}=\stackrel{F}{A}\!{}_{(\rho\beta)},\,\,
\stackrel{F}{A}\!{}^{\rho}{}_{\rho}=0,\\
\stackrel{F}{S}\!{}^{\rho}{}_{\rho\beta}+\stackrel{F}{A}\!{}_{\rho\beta}\ell^{\rho}=0,\,\,
\stackrel{F}{S}\!{}_{\rho\alpha\beta}k^{\rho}=0, \,\,
\stackrel{F}{S}\!{}_{\rho\alpha\beta}k^{\alpha}=0,\,\,
\stackrel{F}{A}\!{}_{\rho\beta}k^{\rho}=0.\hspace{2cm}
\eean
The main condition (\ref{basicC2}) splits then into
\bean
\stackrel{F}{A}\!{}_{\tau}{}^{\rho} \stackrel{G}{A}\!{}_{\rho\beta}&=&0,\\
\stackrel{F}{A}\!{}_{\tau}{}^{\rho} \stackrel{G}{S}\!{}_{\rho\alpha\beta}&=&0,\\
\stackrel{F}{S}\!{}_{(\tau}{}^{\rho}{}_{\nu)}\stackrel{G}{A}\!{}_{\rho\beta}+
\stackrel{F}{A}\!{}_{\tau\nu}\stackrel{G}{A}\!{}_{\rho\beta}\ell^{\rho}&=&0,\\
\stackrel{F}{S}\!{}_{(\tau}{}^{\rho}{}_{\nu)} \stackrel{G}{S}\!{}_{\rho\alpha\beta}+
\stackrel{F}{A}\!{}_{\tau\nu}\stackrel{G}{S}\!{}_{\rho\alpha\beta}\ell^{\rho}&=&0
\eean
for arbitrary 2-forms $F^{\mu\nu}$ and $G^{\mu\nu}$. The first of
these provides, as before,
$$
\stackrel{F}{A}\!{}_{\rho\alpha}=f(F)k_{\rho}k_{\alpha}
$$
but then the third leads to
$$
f(G)\stackrel{F}{A}\!{}_{\tau\nu}=0
$$
which in any case implies 
$$
\stackrel{F}{A}\!{}_{\tau\nu}=0, \,\,\, \forall F\in \Lambda_2\, .
$$
But then one deduces that for arbitrary $F\in \Lambda_2$, 
$\stackrel{F}{\Omega}_{\rho\alpha\beta}=\stackrel{F}{S}\!{}_{\rho\alpha\beta}$ so that
in fact not only $k^{\rho}\nabla_{\rho}C_{\alpha\beta\gamma\delta}=0$
but also
$$
k^{\alpha}\nabla_{\rho}C_{\alpha\beta\gamma\delta}=0
$$
as required.\finn
\begin{lemma}
\label{alg2}
At any point of a Ricci-flat 2-symmetric  Lorentizan manifold $\espaitemps$ where
the Riemann (or equivalently the Weyl)
tensor satisfies (\ref{basicC}) there exists a null vector $k^{\mu}$ such that
\be
\nabla_{\rho}C_{\alpha\beta\gamma\delta}=
4k_{\rho}k_{[\alpha}B_{\beta][\delta}k_{\gamma]}
\label{nablaC}
\ee
where $B_{\beta\mu}$ is a symmetric tensor with the following
properties
\be
k^{\mu}B_{\beta\mu}=0, \,\,\,\, B^{\mu}{}_{\mu}=0.\label{B}
\ee
\end{lemma}
\proof From the previous lemma, we know that,
for {\em arbitrary} tangent vector $\vec v$, the
tensor $\nabla_{\vec v}C_{\alpha\beta\gamma\delta}$ satisfies
$$
\nabla_{\vec v}C_{\alpha\beta\gamma\delta}k^{\alpha}=0
$$
which, with the essential help of (\ref{r2rie}) (in particular using 
$\nabla_{\vec v}C_{\alpha\beta\gamma\delta}\nabla_{\vec v}C^{\alpha\beta\gamma\delta}=0$)
can be shown to imply in our case\footnote{Following the standard
nomenclature in General Relativity \cite{Exact}, a Weyl-like tensor with
this form is said to be of `type III or N', see e.g. \cite{CMPP,CMPP2,PP} and
references therein.} \cite{PP}
$$
\nabla_{\vec v}C_{\alpha\beta\gamma\delta}=4k_{[\alpha}\stackrel{v}{B}_{\beta][\delta}k_{\gamma]}+
2k_{[\alpha}\!\stackrel{v}{D}_{\beta]\gamma\delta}+2k_{[\gamma}\!\stackrel{v}{D}_{\delta]\alpha\beta}
$$
where
\bean
\stackrel{v}{B}_{\beta\mu}=\stackrel{v}{B}_{(\beta\mu)},\,\,
\stackrel{v}{B}_{\beta\mu}k^{\mu}=0,\,\,\,
\stackrel{v}{B}\!{}^{\mu}{}_{\mu}=0, \hspace{25mm}\\
\stackrel{v}{D}_{\beta\lambda\mu}=\stackrel{v}{D}_{\beta[\lambda\mu]},\,\,\,\,
\stackrel{v}{D}_{[\beta\lambda\mu]}=0,\,\,\,\,
k^{\beta}\stackrel{v}{D}_{\beta\lambda\mu}=0,\,\,\,\,
k^{\mu}\stackrel{v}{D}_{\beta\lambda\mu}=0,\,\,\,\,
\stackrel{v}{D}\!{}^{\rho}{}_{\rho\mu}=0
\eean
and, without loss of generality, one can choose another null vector
$\ell_{\mu}$ such that $\stackrel{v}{B}_{\beta\mu}\ell^{\mu}=0$, 
$\ell^{\beta}\stackrel{v}{D}_{\beta\lambda\mu}=0$ and 
$\ell^{\mu}\stackrel{v}{D}_{\beta\lambda\mu}=0$ too. This
leads to
$$
\nabla^{\rho}C_{\alpha\beta\gamma\delta}=
4k_{[\alpha}B^{\rho}{}_{\beta][\delta}k_{\gamma]}+
2k_{[\alpha}D^{\rho}{}_{\beta]\gamma\delta}+2k_{[\gamma}D^{\rho}{}_{\delta]\alpha\beta}
$$
for tensors $B_{\rho\beta\mu}$ and $D_{\rho\beta\lambda\mu}$ with the properties
\bean
B_{\rho\beta\mu}=B_{\rho(\beta\mu)}, \, k^{\rho}B_{\rho\beta\mu}=0,
\,
k^{\beta}B_{\rho\beta\mu}=0, \, B^{\rho}{}_{\rho\mu}=0, \,
B_{\beta}{}^{\rho}{}_{\rho}=0,\, D_{\rho\beta\lambda\mu}=D_{\rho\beta[\lambda\mu]}, \\
D_{\rho[\beta\lambda\mu]}=0,\, k^{\rho}D_{\rho\beta\lambda\mu}=0,\,
k^{\beta}D_{\rho\beta\lambda\mu}=0,\,
k^{\mu}D_{\rho\beta\lambda\mu}=0,\,
D_{\rho}{}^{\sigma}{}_{\sigma\mu}=0,\,
D^{\rho}{}_{\rho\lambda\mu}=0,\, D^{\rho}{}_{\beta\rho\mu}=0
\eean
and, without loss of generality $B_{\rho\beta\mu}\ell^{\mu}=0$,
$D_{\rho\beta\lambda\mu}\ell^{\beta}=0$ and $D_{\rho\beta\lambda\mu}\ell^{\mu}=0$
too. But then, using the last in (\ref{r2rie}) one easily obtains
$$
D^{\rho\beta\lambda\mu}D_{\rho\beta\lambda\mu}=0
$$
which implies necessarily
$$
D_{\rho\beta\lambda\mu}=k_{\rho}A_{\beta\lambda\mu}
$$
for some $A_{\beta\lambda\mu}$ with the necessary properties.
Splitting the other tensor as follows
$$
B_{\rho\beta\mu}=k_{\rho}B_{\beta\mu}+\tilde{B}_{\rho\beta\mu},\,\,\,  
\ell^{\rho}\tilde{B}_{\rho\beta\mu}=0
$$
one can prove, by using  
$\nabla_{[\rho}C_{\alpha\beta]\gamma\delta}=0$, that in fact
$$
A_{\beta\lambda\mu}=2\tilde{B}_{[\lambda\mu]\beta}.
$$
The basic equation (\ref{basicC}) can now be recalculated leading to
$$
\tilde{B}^{\rho}{}_{\beta\delta}B_{\rho\mu}=0, \,\,\,\, 
\tilde{B}^{\rho}{}_{\beta\delta}\left(\tilde{B}_{\rho\mu\nu}-2\tilde{B}_{\mu\nu\rho}\right)=0.
$$
It is an exercise to check that, for a spatial tensor such as
$\tilde{B}_{\rho\mu\nu}$---which belongs to the tensor algebra
constructed over the vector space  $<\vec k,\vec\ell>^{\perp}$--- the
last condition implies
$$
\tilde{B}_{\rho\mu\nu}=0
$$
which leads finally to the result (\ref{nablaC}) with (\ref{B}).\finn

\begin{theorem}
\label{nullorsym}
Let $D\subset \varietat$ be a simply connected domain of an 
$n$-dimensional {\em 2-symmetric} Lorentzian manifold $\espaitemps$. 
Then, if there is no null covariantly constant vector field on $D$, 
$(D,g)$ is in fact locally symmetric.
\end{theorem}
\proof From the theorem (\ref{ric-flat}) one only has to consider the 
case of Ricci-flat manifolds:
\be
R_{\mu\nu}=0, \hspace{2cm} 
R_{\alpha\beta\gamma\delta}=C_{\alpha\beta\gamma\delta}\label{ricciflat}
\ee
so that one can indistinctly use 
the Weyl or the Riemann tensor, while the Ricci tensor has to be set 
equal to zero. From the previous lemma \ref{alg2} it follows that it
must exist a null vector field $k^{\mu}$ and a symmetric tensor field
$B_{\beta\mu}$ on $D$ such that (\ref{nablaC}) and (\ref{B}) hold. It
follows that the tensor
\be
\label{se2}
T_{\alpha\beta\lambda\mu\tau\nu}\equiv 
4\nabla_{\alpha}C_{\lambda\rho\tau\sigma}\nabla_{\beta}C_{\mu}{}^{\rho}{}_{\nu}{}^{\sigma}
=4(B^{\rho\sigma}B_{\rho\sigma})k_{\alpha}k_{\beta}k_{\lambda}k_{\mu}k_{\tau}k_{\nu}
\ee
(this is in fact, in this case, the (super)$^2$-energy
tensor, formula (44) in \cite{S}). As
$T_{\alpha\beta\lambda\mu\tau\nu}$ is parallel by assumption, it
follows that in any region where $B_{\mu\nu}$ is not zero there exists
a parallel null vector field proportional to $k_{\mu}$.\finn

The main results can be summarized in the following central theorem.
\begin{theorem}
\label{classi}
Let $D\subset \varietat$ be a simply connected domain of an 
$n$-dimensional {\em 2-symmetric} Lorentzian manifold $\espaitemps$. 
Then, the line element on $D$ is (possibly a flat extension of) the 
direct product of a certain number of locally symmetric proper 
Riemannian manifolds times either 
\begin{enumerate}
\item a Lorentzian locally symmetric 
spacetime (in which case the whole $(D,g)$ is locally symmetric), or 
\item a Lorentzian manifold with a covariantly constant null vector 
field ---hence belonging to the general Brinkmann's class presented below in
(\ref{brink}).
\end{enumerate}
\end{theorem}
\proof
It follows from Lemma \ref{redornull} and the 
previous theorems \ref{ric-flat} and \ref{nullorsym}.\finn

\noindent
{\bf Remarks:}
\begin{enumerate}
\item Of course, the number of proper Riemannian symmetric manifolds 
can be zero, so that the whole 2-symmetric spacetime, if not locally 
symmetric, is given just by a line-element of the form (\ref{brink}).
\item Although mentioned explicitly for the sake of clarity, it is 
obvious that the block added in any flat extension can also be 
considered as a particular case of a locally symmetric part building 
up the whole space.
\item This theorem provides a full characterization of the 
2-symmetric spaces using the classical results on the symmetric ones: 
their original classification (for the semisimple case) was given in 
\cite{C1}, see also \cite{Ber,H}, and the general problem was solved for Lorentzian 
signature in \cite{CW}. Combining these results with those for proper 
Riemannian metrics \cite{C,C1,H}, a complete classification is 
achieved.
\end{enumerate}

The most general 
local line-element for a Lorentzian manifold with covariantly constant
null vector field was discovered by Brinkmann 
\cite{Br} by studying the Einstein spaces which can be mapped 
conformally to each other. In appropriate local coordinates 
$$
\{x^0,x^1,x^i\}=\{u,v,x^i\},\,\,\,\, (i,j,k,\dots =2,\dots ,n-1)
$$
the line-element reads (see also \cite{Z})
\be
ds^2=-2du(dv+Hdu+W_i dx^i)+g_{ij}dx^idx^j \label{brink}
\ee
where the functions $H$, $W_i$ and $g_{ij}=g_{ji}$ are independent of 
$v$, otherwise arbitrary, and the parallel null vector field is given 
by
\be
k_{\mu}dx^{\mu}=-du, \hspace{1cm} k^{\mu}\partial_{\mu}=\partial_v \, .
\ee
Given that all 2-symmetric non-symmetric Lorentzian manifolds contain 
a covariantly constant null vector field according to theorems
\ref{nullorsym} and \ref{classi}, and that the explicit local form
(\ref{brink}) is known, 
it is a simple matter of calculation to identify which manifolds among (\ref{brink}) 
are actually 2-symmetric. Using theorem \ref{classi} and its remarks, 
this provides ---by direct product with proper Riemannian symmetric 
manifolds if adequate--- all possible {\em non-symmetric} 2-symmetric 
spacetimes. This is the subject of the second paper \cite{S2}.

A particularly simple example of these 2-symmetric spacetimes is
provided by the ``plane waves'' (see e.g. \cite{EK,Exact}) for which\footnote{Actually, in $n=4$, all 2-recurrent spacetimes were 
given in \cite{MTh}, see also \cite{Th,Th1,Th2} where they were proved to be a subset of the Brinkmann class (\ref{brink}), and the conformally 2-recurrent spacetimes belonging to (\ref{brink}) and satisfying $W_i=0$, $g_{ij}=\delta_{ij}$, were all found in \cite{CS}.} $g_{ij}=\delta_{ij}$, $W_i=0$ and the remaining function is a quadratic
expression on the coordinates $x^i$
\be
H=a_{ij}(u)x^ix^j \, .\label{pw}
\ee
The functions $a_{ij}(u)$ are arbitrary for general plane waves, but in our case
the 2-symmetry is easily seen to imply then that
$$
\frac{d^2a_{ij}}{du^2}=0,\,\,\,\, \Longrightarrow \,\,\,\,\, a_{ij}(u)=c_{ij}u+b_{ij}
$$
for some constants $c_{ij},b_{ij}$. These spaces have been actually
identified recently in \cite{PCM}. Note that one can add any number of
proper Riemannian locally symmetric spaces to this one keeping the
whole 2-symmetry.

\section{Future work and open problems}
Arguably, the decisive step to characterize all $k$-symmetric
spacetimes has been achieved with theorems \ref{ric-flat},
\ref{nullorsym} and \ref{classi}. The same or even more modest
techniques ---for instance, using repeatedly the Ricci identity---
will allow for a complete solution to the problem with
arbitrary $k$ now that the 2nd order stage has been solved. As a
matter of fact, it seems plausible that the following conjecture will
hold
\begin{conj}
\label{conjec}
Let $D\subset \varietat$ be a simply connected domain of an 
$n$-dimensional $k$-symmetric Lorentzian manifold $\espaitemps$. 
Then, if there is no null covariantly constant vector field on $D$, 
$(D,g)$ is in fact locally symmetric.
\end{conj}
If this is true, then a result similar to theorem \ref{classi} will be
valid and the $k$-symmetric spacetimes will be just the direct sum of
lower-order symmetric ones plus a particular $k$-symmetric one of the Brinkmann
class (\ref{brink}). Actually, there are obvious
$k$-symmetric spaces within the plane waves: just use the same form
(\ref{pw}) for the function $H$ but let the functions $a_{ij}(u)$
be polynomials of degree up to $k-1$, see also \cite{PCM}. 

As a matter of fact, it might be the case that not only
$k$-symmetric, but also all $k$-recurrent spacetimes, have the same
type of hierarchy, see also \cite{CS}. I would like to put forward the second conjecture:
\begin{conj}
\label{conjec2}
Let $D\subset \varietat$ be a simply connected domain of an 
$n$-dimensional $k$-recurrent Lorentzian manifold $\espaitemps$. 
Then, if there is no null recurrent vector field on $D$, 
$(D,g)$ is in fact recurrent.
\end{conj}
In this sense, it should be interesting to know if there is any
characterization of the spacetimes which belong to the Brinkmann
class ---that is, they have a parallel null vector field---, or to the more general classes with a recurrent null vector field, see e.g. \cite{Leis2}, but are
neither recurrent nor symmetric of any order, if they exist.

Of course, all the above can be generalized, if desired, to
conformally $k$-symmetric, or $k$-recurrent,  and Ricci $k$-symmetric
and $k$-recurrent.

There are several other routes open and worth exploring. To start
with, this paper has concentrated on Lorentzian signature, but the
result does not apply, in principle, to other
signatures. Mathematically, the resolution of the second-order
symmetric spaces will not be complete until these other cases are
fully solved. 

Another interesting problem is that of semi-symmetric spacetimes (and
with other signatures). This problem, as well as the conformally
semi-symmetric one, can be solved in 4-dimensional Lorentizan manifolds by using spinors \cite{PR}, see the classification in \cite{HV}.
However, as happened with the case treated herein, the
higher dimensional cases need further exploration.

Along the paper some relationships with holonomy, parallel tensors and
vectors, or scalar invariants
have appeared. In this sense, it would be interesting to know the list
of all spacetimes with {\em constant} curvature scalar
invariants \cite{CHP}. Similarly, the implications of the existence of parallel
tensors with rank higher than 2 would be of enormous help in this and
related studies---apart from interesting on its own
right---. 

Finally, the algebraic solution of equations such as
(\ref{RR}), or more interestingly,
such as (\ref{RicR}), (\ref{RicC}) and (\ref{CC}) in general
dimension $n$ would be of great help.

\section*{Acknowledgements} 
Comments and useful information from M. S\'anchez, A.J. Di Scala, R. Vera, S.B. Edgar and a referee are gratefully 
acknowledged. Inspiration for some of the results in section \ref{generic}
was gained from manuscript papers ---which I lost and have been unable to find
anywhere--- kindly provided by C.B. Collinson some years ago.
Financial support from the Wenner-Gren Foundations, 
Sweden, is gratefully acknowledged. I thank the Applied Mathematics 
Department at Link\"oping University, where this work was mainly carried 
out, for hospitality. I am also grateful to the School of Mathematical
Sciences, Queen Mary (London, UK) where this manuscript was
finalized in a recent short visit.
Supported by grants FIS2004-01626 (MEC) and GIU06/37 (UPV/EHU).

\end{document}